\documentclass[twoside]{article}
\usepackage{amsfonts,amsmath,amssymb,amscd}
\usepackage{latexsym}
\usepackage{euscript}
\usepackage{color}
\usepackage[T2A]{fontenc}
\usepackage[cp1251]{inputenc}
\usepackage[english,russian]{babel}
\usepackage{multicol}
\usepackage{graphicx}

\newcommand{\n}[1]{\refstepcounter{equation}\label{#1} \eqno{(\arabic{section}.\arabic{equation})}}

\makeatletter
\renewcommand{\@oddfoot}{}%
\renewcommand{\@evenfoot}{}%
\renewcommand{\@oddhead}{%
\raisebox{0pt}[\headheight][0pt]{\vbox{\hbox to
\textwidth{\strut%
\hfill\rightmark\hfill\thepage%
}\hrule}}%
}
\renewcommand{\@evenhead}{%
\raisebox{0pt}[\headheight][0pt]{\vbox{\hbox to \textwidth{\strut%
\large\thepage\hfill\leftmark\hfill%
}\hrule}}%
}

\begin{document}
\begin{Large}

\newtheorem{theorem}{ \sc Теорема}
\newtheorem{lemma}{\sc Лемма}
\newtheorem{statement}{\sc Утверждение}
\newtheorem{proposition}{\sc Предложение}
\newtheorem{corollary}{\sc Следствие}
\newtheorem{remark}{\sc Замечание}
\newtheorem{definition}{\sc Определение}
\newtheorem{example}{\sc Пример}

\@addtoreset{equation}{section}
\renewcommand{\@evenhead}{\large\thepage\hfill}
\renewcommand{\@oddhead}{\large\hfill\thepage}
\makeatother

\newcommand{\sect}[1]{\refstepcounter{section}
\begin{center}
\LARGE \bf\vspace*{3mm} \S \thesection. \LARGE\bf #1
\end{center}
\vskip3mm}


\renewcommand{\abstract}[1]{\begin{center}
\begin{minipage}[t]{150mm}
\small #1
\end{minipage}
\end{center}}

\newcommand{\head}[2]{\markboth{\large #1}{\large #2}}

\renewcommand{\refname}{%
                          \vspace*{5mm}
                          \centerline {\bf Список литературы}
                          \vspace*{-10mm}
                       }

\def\Re{\operatorname{Re}}
\def\Im{\operatorname{Im}}
\def\const{\mathrm{const}}
\def\RR{\mathbb R}
\def\CC{\mathbb C}
\def\NN{\mathbb N}
\def\RS{\mathfrak R}
\def\bz{\mathbf z}
\def\sH{\mathscr H}
\def\HH{\mathscr H}
\def\mro{\widehat\rho}
\def\FR{{\cal F}({\Bbb R})}
\def\FRn{{\cal F}({\Bbb R}^n)}
\def\Rn{{\Bbb R}^n}
\def\KR{{\cal K}({\Bbb R})}
\def\KRn{{\cal K}({\Bbb R}^n)}
\def\MR{M({\Bbb R}^n)}
\def\W{\cal W}
\head{\sc А.\,Ф.~Гришин,  М.\,В.~Скорик}%
{\sc Некоторые свойства интеграла Фурье} \thispagestyle{empty}
\begin{center}
{\LARGE\bf A.\,F.~Grishin, M.\,V.~Skoryk }
 \vskip3mm
{\LARGE\bf Some properties of Fourier integrals}
\end{center}

\begin{abstract}
{\large {\bf MSC subject classification.} 42A38

Let ${\cal F}({\Bbb R}^n)$ be the algebra of Fourier transforms of
functions from $L_1(\RR^n),$ $\KRn$ be the algebra of Fourier
transforms of bounded complex Borel measures in $\Rn$  and ${\cal
W}$ be Wiener algebra of continuous $2\pi$-periodic functions with
absolutely convergent Fourier series. New properties of functions
from these algebras are obtained.

Some conditions which determine membership of $f$ in $\FR$ are
given. For many elementary functions $f$ the problem $f\in\FR$ can
be resolved easily using these conditions. We prove that the
Hilbert operator is a bijective isometric operator in the Banach
spaces $\cal W_0,$ $\FR,$ $\KR-A_1$ ($A_1$ is the one-dimension
space of constant functions). We also consider the classes $M_k,$
which are similar to the Bochner classes $F_k,$ and obtain
integral representation of the Carleman transform of measures of
$M_k$ by integrals of the form
$\displaystyle\int\limits_{-\infty}^\infty\frac{u(\lambda)}{(\lambda+z)^{k+1}}\;d\lambda.$

References: 22 units.

{\bf Key words:} Wiener algebra, Krein algebra,  Hilbert
transformation, Carleman trans\-formation, Povzner formula.

{\bf Comments:} 36 pages, in Russian}
\end{abstract}

\begin{abstract}
{\large {\bf УДК} 517.443

Пусть ${\cal F}({\Bbb R}^n)$ --- алгебра функций, являющихся
преобразованиями Фурье функций из $L_1(\RR^n)$, а $\KRn$ ---
алгебра функций, являющихся преобразованиями Фурье конечных
комплексных борелевских мер в $\Rn$, ${\cal W}$ --- винеровская
алгебра непрерывных $2\pi$-периодических функций с абсолютно
сходящимся рядом Фурье. В статье приводятся несколько новых
свойств функций из этих алгебр.

Доказываются несколько различных условий принадлежности функций
$\varphi$ алгебре $\FR.$ С помощью этих условий для многих
элементарных функций задача о принадлежности их к алгебре $\FR$
решается достаточно просто. Доказывается, что оператор Гильберта
--- это биективный изометрический оператор в банаховых
пространствах $\cal W_0,$ $\FR,$ $\KR-A_1$($A_1$ --- одномерное
пространство постоянных функций). Это является основным
результатом работы. Кроме того, дается интегральное представление
преобразований Карлемана мер из классов $M_k,$ аналогичных классам
функций $F_k,$ введённых Бохнером, интегралами вида
$\displaystyle\int\limits_{-\infty}^\infty\frac{u(\lambda)}{(\lambda+z)^{k+1}}\;d\lambda.$

Библиография: 22 названия.

{\bf Ключевые слова:} алгебра Винера, алгебра Крейна,
преобразование Гильберта, преобразование Карлемана, формула
Повзнера.}

\end{abstract}

\sect{Вступление}\label{vstuplenie2}

Литература по теории преобразования Фурье многими характеризуется
как необозримая. Мы приведём некоторые из результатов этой теории.

{\bf 1.} Если $f\in L_1(-\infty,\;\infty),$ то функция
$\displaystyle\widehat f
(x)=\frac{1}{\sqrt{2\pi}}\int\limits_{-\infty}^{+\infty}f(t)e^{-ixt}dt$
непрерывна на оси $(-\infty, +\infty)$. \par {\bf 2.} Если $f\in
L_1(-\infty,\;\infty),$ то $\lim\limits_{x\to\pm\infty}\widehat f
(x)=0$.

{\bf 3.} Если преобразование Фурье $\widehat f$ функции $f\in
L_1({\Bbb R}^n)$ определить формулой $$\widehat
f(x)=\left(\frac{1}{\sqrt{2\pi}}\right)^n\int\limits_{{\Bbb R}^n}
f(t)e^{-i(x, t)}dt,$$ то для преобразования Фурье свертки
$$f(t)=(f_1*f_2)(t)=\int\limits_{{\Bbb R}^n} f_1(t-u)f_2(u)du
$$
справедлива формула $\widehat
f(x)=\left(\sqrt{2\pi}\right)^n\widehat f_1(x)\widehat f_2(x).$
Это важное свойство преобразований Фурье. Из него следует, что
множество функций $\FRn$ является алгеброй относительно поточечных
умножения и сложения.

{\bf 4.} Если $M({\Bbb R}^n)$ --- алгебра конечных комплексных
борелевских мер в $\Rn,$ $\mu\in\MR,$ а $\widehat\mu(x)$
--- преобразование Фурье меры $\mu,$ то $\widehat\mu(x)$ ---
ограниченная равномерно непрерывная функция в $\Rn.$ Множество
функций вида $\widehat\mu(x)$ образует алгебру относительно
поточечных сложения и умножения. Мы будем обозначать эту алгебру
через $\KRn$ и называть алгеброй Крейна.

{\bf 5.} Теорема Бохнера утверждает, что класс непрерывных
положительно определённых функций в $\Rn$ совпадает с классом
преобразований Фурье положительных мер из $\MR.$ Отсюда, в
частности, следует, что алгебра Крейна $\KRn$ совпадает с
комплексной линейной оболочкой множества непрерывных положительно
определённых функций в $\Rn.$

{\bf 6.} Если $\psi(x)$ --- чётная непрерыная функция, убывающая и
выпуклая на полуоси $[0, \infty)$, причём $\psi(x)\to0$
($x\to+\infty$), то $\psi\in\FR.$ Это следствие из теоремы~124,
\cite{Skoryk-Titchmarsh}.

{\bf 7.} Обозначим через $\W$ алгебру Винера, состоящую из
непрерывных  $2\pi$-пе\-ри\-о\-ди\-чес\-ких функций, которые
разлагаются в абсолютно сходящийся ряд Фурье. Известно, что
алгебры $\W,$ $\FR,$ $\KR$ локально совпадают. Это означает
следующее. Пусть $f_1$
--- функция, принадлежащая одной из названных выше трёх алгебр,
$x\in{\Bbb R}.$ Тогда существует число $\delta>0$ и функции $f_2$
и $f_3,$ принадлежащие двум другим алгебрам, такие, что при
$t\in(x-\delta,\;x+\delta)$ выполняется равенство
$f_1(t)=f_2(t)=f_3(t).$ Впервые этот факт заметил
Винер~(\cite{Skoryk-Wiener}). Однако, приведенная формулировка в
книге Винера отсутствует. Доказательство приведенного утверждения
и дальнейшие ссылки можно найти в \cite[глава 2, пункт
4]{Skoryk-Kahan}.

{\bf 8.} Теорема Винера утверждает, что если $2\pi$-периодическая
функция $\varphi$ локально принадлежит $\W,$ то $\varphi\in\W.$
Доказательство теоремы Винера можно найти в \cite[глава 2, пункт
4]{Skoryk-Kahan}.

{\bf 9.} Справедлив следующий результат Стечкина. Пусть $f\in
L_2(0,\; 2\pi),$
$$e_n(f)=\inf\|f(t)-\sum_{k=1}^n c_k e^{i\lambda_k t}\|_{L_2(0, 2\pi)},$$ $\lambda_k$ --- вещественные числа,
$c_k$ --- комплексные числа. Тогда для того, чтобы функция $f$
была эквивалентна функции из алгебры $\W$ необходимо и достаточно,
чтобы $$\displaystyle\sum_{n=1}^\infty
\frac{e_n(f)}{\sqrt{n}}<\infty.$$

Доказательство теоремы Стечкина можно найти в \cite[глава 2, пункт
3]{Skoryk-Kahan}.

{\bf 10.} Сформулируем ещё достаточные условия Бернштейна и
Зигмунда принадлежности функции алгебре $\W.$ Пусть $\Phi$
--- класс возрастающих на полуоси $[0,\infty)$ функций $\varphi,$
удовлетворяющих условиям: $\varphi(0)=0,$
$\varphi(x+y)\le\varphi(x)+\varphi(y).$

Пусть $f$ --- непрерывная $2\pi$-периодическая функция,
$\omega_f(\delta)$ --- её модуль непрерывности, причем
$\omega_f(\delta)\le\varphi(\delta),\quad\varphi\in\Phi.$ Тогда
если
$$\int\limits_0^1\frac{\varphi(x)}{x^\frac{3}{2}}dx\;<\;\infty, $$
то $f\in\W.$

Пусть $f$ --- непрерывная $2\pi$-периодическая функция, имеющая
ограниченную вариацию на сегменте $[0,\; 2\pi].$ Пусть выполняется
неравенство
$\omega_f(\delta)\le\varphi(\delta),\quad\varphi\in\Phi.$ Тогда из
сходимости интеграла $$\int\limits_0^1
\frac{\sqrt{\varphi(x)}}{x}dx$$ следует, что $f\in\W.$

Доказательство теорем Бернштейна и Зигмунда также можно найти в
\cite[глава 2, пункт 6]{Skoryk-Kahan}.

{\bf 11.} Для того, чтобы функция $F\in L_q(-\infty,\;\infty)$
была преобразованием Фурье некоторой функции из
$L_p(-\infty,\;\infty),$ $p\in(1,2],$ $\displaystyle
\frac{1}{p}+\frac{1}{q}=1,$ необходимо  и достаточно,  чтобы для
функции
$$\varphi(x)=\int\limits_{-\infty}^\infty F(t)\frac{e^{itx}-1}{it}\;dt$$
выполнялись неравенства
$$\sum\limits_{k=1}^{n-1}\frac{|\varphi(x_{k+1})-\varphi(x_k)|^p}{(x_{k+1}-x_k)^{p-1}}\le M(F)$$
для любых $x_k$ таких, что $-\infty<x_1<x_2<\ldots<x_n<\infty.$

Это теорема 1 из \cite{Skoryk-Berry}.

{\bf 12.} Для того, чтобы функция $F\in L_\infty(-\infty,\infty)$
принадлежала алгебре $\FR,$ необходимо и
 достаточно, чтобы выполнялись
следующие два условия:

1) функционал $$\lambda_F(\widehat g)=\int\limits_{-\infty}^\infty
F(x)g(x)dx $$ был непрерывным линейным функционалом на линейном
многообразии пространства $L_\infty(-\infty,\;\infty),$ состоящем
из функций алгебры $\FR,$

2) $$\lim\limits_{n\to\infty}\int\limits_{-\infty}^\infty
F(x)g_n(x)dx=0$$ для любой последовательности функций $g_n\in
L_1(-\infty,\infty)$ такой, что $\|\widehat g_n\|_1{\to}0$
$(n\to\infty),$ $\|\widehat g_n\|_\infty\le 1.$

Это переформулировка теоремы 2 из \cite{Skoryk-Berry}.

{\bf 13.} Для того, чтобы непрерывная ограниченная функция $f(t)$
на оси $(-\infty, \;\infty)$ принадлежала алгебре $\KR,$
необходимо и достаточно, чтобы выполнялось неравенство
$$\displaystyle\sup\limits_{n\ge1}\int\limits_{-\infty}^\infty\left|\int\limits_{-\infty}^\infty
f(t)\left(\frac{\sin\frac{t}{n}}{\frac{t}{n}}\right)^2e^{-it\lambda}\;dt\right|d\lambda<\infty.$$

Это с точностью до терминологии  из \cite[теорема
3]{Skoryk-Yosida}.

{\bf 14.} Рассмотрим ряд $$\varphi_\alpha(x)=\sum_{n=1}^\infty
e^{in\; {\rm ln}\; n}\frac{e^{inx}}{n^{\frac{1}{2}+\alpha}},
\quad\alpha\in(0,1).$$

Написанный ряд сходится равномерно на всей оси $(-\infty,\infty),$
а функция $\varphi_\alpha(x)$ удовлетворяет условию Гёльдера
порядка $\alpha.$ Это утверждение приведено в книге Зигмунда
\cite[глава 5, раздел 4]{Skoryk-Zigmund-1}. В той же книге в
примечаниях имеются ссылки на первоисточники.

Отметим ещё, что в книге Эдвардса \cite[глава 10, пункт
10.6]{Skoryk-Edwards} приведены различные результаты, касающиеся
алгебры $\W$ и приведено большое количество ссылок на работы по
этой тематике.

В параграфе \ref{suff} даются различные достаточные условия
принадлежности функций алгебре $\FR.$ Например, условие
выпуклости, которое присутствует в сформулированном выше
утверждении 6, заменяется на условие кусочной выпуклости. Причем
направления выпуклости на различных участках не связаны между
собой. Из приведенных в параграфе~2 условий, в частности, следует
принадлежность функции $\displaystyle\frac{1}{\ln_k (\alpha_k
+|x|)},$ где $\ln_kx=\ln\ldots\ln x,$ $\alpha _k$
--- достаточно большое число, алгебре $\FR.$ Здесь $k\ge 1$ ---
любое целое число.

Отметим некоторые результаты из параграфа~\ref{razdel_KR}.

Из сформулированного утверждения 6 следует, что в алгебре $\FR$
есть функция $\varphi(x),$ которая совпадает с функцией
$\displaystyle\frac{1}{\ln_kx}$ в некоторой окрестности $+\infty.$
Можно поставить вопрос: существует ли в алгебре $\FR$ функция
$\varphi$ такая, что
$\displaystyle\varphi(x)\sim\frac{2}{\ln_2x}\;(x\to+\infty),\;\;\varphi(x)\sim
\frac{3}{\ln_3|x|}\;(x\to-\infty)?$ Ответ отрицательный.
Оказывается, что функции $\varphi$ из $\FR,$ и даже из $\KR,$
''почти'' четные. Точнее, мы доказываем (теорема
\ref{Skoryk-Theorem-6}), что если $\varphi\in\KR$, а $z$ ---
невещественное число, то существует $\displaystyle
V.P.\int\limits_{-\infty}^{+\infty}\frac{\varphi(\lambda)}{\lambda
+ z}d\lambda.
$
Это довольно существенное ограничение на поведение функций
$\varphi$  в окрестности бесконечности, если учитывать, что
функции из $\KR$ не обязательно бесконечно малые на бесконечности.

Точка $\infty$ не является исключением. Существуют ограничения и
на локальное поведение функций $\varphi$ из алгебры $\KR$.

Мы доказываем (теоремы~\ref{Skoryk-Theorem-9},
\ref{Skoryk-Theorem-H4}), что функция
$$\psi(x)=V.P.\frac{1}{\pi}\int\limits_{-\infty}^{+\infty}\frac{\varphi(\lambda)}{x-\lambda}d\lambda,$$
где интеграл понимается, как
$$\lim\limits_{{N\to+\infty}\atop{\varepsilon\to
+0}}\frac{1}{\pi}\left(\int\limits_{-N}^{x-\varepsilon}\frac{\varphi(\lambda)}{x-\lambda}d\lambda
+
\int\limits_{x+\varepsilon}^{N}\frac{\varphi(\lambda)}{x-\lambda}d\lambda\right),$$
корректно определена для всех $x\in(-\infty, +\infty)$ и
принадлежит алгебре $\KR$($\FR,$ если $\varphi\in\FR$).
Аналогичный результат (теорема~\ref{Skoryk-Theorem-H8}) справедлив
и для алгебры $\W.$ Мы доказываем, что если $\varphi\in\W,$ то
$$V.P.\frac{1}{\pi}\int\limits_{\Bbb T}\frac{\varphi(\zeta)}{\zeta-z}d\zeta\;\in\W,$$
где ${\Bbb T}=\left\{\zeta:\;\;|\zeta|=1\right\},\;\;z\in{\Bbb
T}.$

Остановимся на результатах параграфа \ref{4}. Пусть $k\ge 0$ ---
целое число. Класс Бохнера $F_k$ состоит из функций $f(t)$
измеримых на вещественной оси и таких, что выполняется условие
$\displaystyle\frac{f(t)}{1+|t|^k}\in L_1(-\infty, +\infty)$. Для
$f\in F_k$ определяется преобразование Карлемана $F(z) = (F_+(z),
F_-(z)),$ где $$F_+(z) = \int\limits_0^{+\infty}f(t)e^{itz}dt,
\quad\Im z > 0,\quad F_-(z) = -\int\limits_{-\infty}^0
f(t)e^{itz}dt,\quad \Im z < 0.$$

Функция $F(z)$ является кусочно аналитической функцией в
комплексной плоскости, разрезанной по вещественной оси.

Мы вводим классы $M_k$ мер, которые аналогичны классам $F_k.$ Мера
$\mu$ принадлежит классу $M_k,$ если мера $\mu_1,$ $\displaystyle
d\mu_1(x) = \frac{d\mu(x)}{1+|x|^k},$ принадлежит классу $M({\Bbb
R}).$ Для мер из класса $M_k$ также определяется преобразование
Карлемана
$$F_+(z) = \int\limits_0^{+\infty \;\;_\prime} e^{itz}
d\mu(t),\quad \Im z >0,\quad F_-(z)
=-\int\limits_{-\infty}^{\;\;0\;\;\;_\prime} e^{itz}d\mu(t),\quad
\Im z <0, $$ где штрихи над знаками интегралов означают, что
интегралы по полуосям $[0, +\infty)$ и $(-\infty, 0]$ берутся не
по мере $\mu,$ а по мере $\mu-\frac{1}{2}\mu(\{0\})\delta,$ где
$\delta$ --- мера Дирака.

В параграфе~\ref{4} даются представления интегралами вида
$\displaystyle\int\limits_{-\infty}^{+\infty}\frac{u(\lambda)}{(\lambda+z)^{k+1}}\,d\lambda$
преобразований Карлемана $F(z)$ мер из классов $M_k$. Доказываемая
формула является обобщением формулы Повзнера
\cite{Skoryk-Povzner}.

Преобразование Карлемана --- важный объект в гармоническом
анализе. Он был предметом исследований для многих математиков.
Карлеман ввёл свое преобразование в \cite{Skoryk-Carleman}. О
развитии идей Карлемана можно прочитать в \cite{Skoryk-Gurariy}
--- \cite{Skoryk-Domar}, где
так же имеются многочисленные ссылки.

\sect{Достаточные условия принадлежности функции алгебре
$\FR$}\label{suff}

В утверждении 6, сформулированном во вступлении, приводится
достаточное условие, найденное Титчмаршем, которое гарантирует
принадлежность заданной функции алгебре $\FR.$ Далее приводятся
более слабые ограничения на функцию $R(x),$ гарантирующие
включение $R\in\FR.$ Полученные результаты позволяют для многих
элементарных функций достаточно просто решать задачу о
принадлежности этих функций алгебре $\FR.$ Конечно, для
элементарных функций $R$ такая задача решается и более простыми
методами, основанными на использовании таблиц преобразований Фурье
и соответствующих асимптотических формул. Наши результаты
позволяют другим методом доказывать соотношение $R\in\FR$ для
таких функций.

Мы начнем со следующего утверждения.

\begin{lemma}\label{Skoryk-Lemma-f1}\hskip-2mm.
Пусть $h\in\ L_1(0,\;\infty),\;\;\alpha\in[0,\;\pi],\;k\ge0$ ---
целое число,
$$H_1(t)=\int\limits_{\textstyle \frac{k\pi}{t}}^{\textstyle \frac{k\pi+\alpha}{t}} h(x)\sin xt \;dx,
\qquad
H_2(t)=\int\limits_{\textstyle\frac{(k+\frac{1}{2})\pi}{t}}^{\textstyle\frac{(k+\frac{1}{2})\pi+\alpha}{t}}
h(x)\cos xt \;dx.$$ Тогда справедливы неравенства
$$\int\limits_0^\infty\frac{|H_1(t)|}{t}\;dt\le\int\limits_0^\alpha \frac{\sin v}{v+k\pi}dv
\int\limits_0^\infty|h(x)|dx,$$
$$\int\limits_0^\infty\frac{|H_2(t)|}{t}\;dt\le\int\limits_0^\alpha \frac{\sin v}{v+(k+\frac{1}{2})\pi}dv
\int\limits_0^\infty|h(x)|dx.$$
\end{lemma}

{\sc Доказательство.} Имеем
$$H_1(t)=\frac{1}{t}\int\limits_{k\pi}^{k\pi+\alpha}h\left(\frac{u}{t}\right)\sin u\;du=\frac{(-1)^k}{t}
\int\limits_0^\alpha h\left(\frac{v+k\pi}{t}\right)\sin v\;dv.$$
Поэтому
$$\int\limits_0^\infty\frac{|H_1(t)|}{t}\;dt\leq\int\limits_0^\alpha\sin v\int\limits
_0^\infty \frac{1}{t^2}|h\left(\frac{v+k\pi}{t}\right)|\;dt\;dv
=\int\limits_0^\alpha\frac{\sin v}{v+k\pi}\int\limits_0^\infty
|h(\tau)|\;d\tau\;dv.$$ Аналогично оценивается интеграл с функцией
$H_2(t).$

Лемма доказана.

Отметим, что приведенные рассуждения не позволяют оценивать более
общие функции $H_1(t)$ и $H_2(t),$ которые отличаются от
рассмотренных в лемме функций $H_1(t)$ и $H_2(t)$ тем, что
целочисленная постоянная $k\ge0$ заменяется на целочисленную
функцию $k(t)\ge0.$

Следующие теоремы посвящены преобразованиям Фурье монотонных
функций. Результаты на эту тему и ссылки на первоисточники, в
частности, на работы Принсгейма можно найти, например, в книге
\cite[раздел 1.10]{Skoryk-Titchmarsh}.

\begin{theorem}\label{Skoryk-Theorem-S1} \hskip-2mm.
Пусть $\varphi(x)$ --- убывающая на полуоси $(0,\;\infty)$ функция
со сходящимся интегралом $\displaystyle\int\limits_0^\infty
\varphi(x)\;dx,$ $\displaystyle
g_m(t)=\int\limits_{\textstyle\frac{m\pi}{t}}^\infty
\varphi(x)\sin xt\;dx.$

Тогда $(-1)^m g_m(t)$ есть положительная на полуоси $(0,\infty)$
функция и выполняется равенство
$$\int\limits_0^\infty\frac{g_m(t)}{t}\;dt=\int\limits_{m\pi}^{\infty}\frac{\sin u}{u}\;du
\int\limits_0^\infty\varphi(x)\;dx.$$

В частности,
$$
\int\limits_0^\infty\frac{g(t)}{t}\;dt=\frac{\pi}{2}\int\limits_0^\infty\varphi(x)\;dx,\quad(g(t)=g_0(t)).\n{formula-10081}
$$
\end{theorem}

{\sc Доказательство.} Будем считать, что $m=0.$ Для других $m$
проходит аналогичное доказательство. Из условий теоремы легко
следует, что $\varphi$
--- положительная бесконечно малая на бесконечности функция.
Кроме того,  $\varphi\in L_1(0,\;\infty).$ Имеем
$$g(t)=\frac{1}{t}\int\limits_0^\infty\varphi\left(\frac{u}{t}\right)\sin u\;du=
\frac{1}{t}\sum\limits_{k=1}^{\infty}\int\limits_{(k-1)\pi}^{k\pi}\varphi\left(\frac{u}{t}\right)
\sin u\;du=$$

$$
=\frac{1}{t}\sum\limits_{k=1}^\infty(-1)^{k-1}\int\limits_0^\pi\varphi\left(
\frac{v+(k-1)\pi}{t}\right)\sin v\; dv.\n{formula-10082}
$$
Обозначим
$$a_k=\int\limits_0^\pi\varphi\left(
\frac{v+(k-1)\pi}{t}\right)\sin v\; dv.$$

Последовательность $a_k$ --- это убывающая бесконечно малая
последовательность. По теореме Лейбница функция $g(t)$
положительна на полуоси $(0,\;\infty)$ и для остатка ряда
(\ref{formula-10082}) справедлива оценка
$$\left|R_n(t)\right|\le\frac{1}{t}\int\limits_0^\pi\varphi\left(\frac{v+n\pi}{t}\right)\sin v\;dv.$$

Тогда, повторяя соответствующую выкладку из доказательства
леммы~\ref{Skoryk-Lemma-f1}, получим
$$\int\limits_0^\infty\frac{|R_n(t)|}{t}\;dt\le\int\limits_0^\pi\frac{\sin v}{v+n\pi}\;dv
\int\limits_0^\infty\varphi(x)\;dx.$$ Это неравенство оправдывает
перемену порядка суммирования и интегрирования в нижеследующей
выкладке
$$\int\limits_0^\infty\frac{g(t)}{t}\;dt=\sum\limits_{k=1}^\infty(-1)^{k-1}\int
\limits_0^\pi\sin v\int\limits_0^\infty
\frac{1}{t^2}\;\varphi\left(\frac{v+(k-1)\pi}{t}\right)\;dt\;dv=$$
$$=\int\limits_0^\infty\varphi(x)\;dx\sum\limits_{k=1}^\infty(-1)^{k-1}\int
\limits_0^\pi\frac{\sin
v}{v+(k-1)\pi}\;dv=\int\limits_0^\infty\varphi(x)\;dx\int\limits_0^\infty\frac{\sin
u}{u}\;du.$$

Тем самым теорема доказана.

\begin{remark} \hskip-2mm.
 Выкладки, применённые при доказательстве формулы
 (\ref{formula-10081}) восходят к Лобачевскому
 \cite{Skoryk-Lobachevskiy-1} (см. также
 \cite{Skoryk-Lobachevskiy-2}), который вычислял интеграл $\displaystyle\int\limits
 _0^\infty f(x)\frac{\sin x}{x}\;dx.$ Эти же выкладки применял
 Титчмарш. Наше доказательство положительности $g(t)$ совпадает с
 доказательством теоремы 123 из книги Титчмарша. Титчмарш
 накладывал более слабое ограничение
 $\displaystyle\int\limits_0^1\varphi(x)dx<\infty$ на функцию
 $\varphi(x).$ Однако, в этом случае нельзя писать равенство
 (\ref{formula-10081}).   Таким образом, рассуждения Титчмарша
 показывают, что положительность функции $g(t)$ следует из более
 слабого ограничения на функцию $\varphi(x),$ чем то, которое
 требуется в теореме. Однако в дальнейшем для нас важно, чтобы
 сходился интеграл
 $\displaystyle\int\limits_0^\infty\frac{g(t)}{t}\;dt.$
 Доказательство теоремы~\ref{Skoryk-Theorem-S1} --- это лишь
 слегка измененные рассуждения Титчмарша.

\end{remark}

 Далее символом $f'_+(x)$ обозначается правая
 производная функции $f$ в точке~$x.$

\begin{theorem}\label{Skoryk-Theorem-S2} \hskip-2mm.
Пусть $f(x)$ --- убывающая выпуклая бесконечно малая в
бесконечности функция на полуоси $[0,\;\infty)$ и пусть
$$F_c(t)=\sqrt\frac{2}{\pi}\int\limits_0^\infty f(x)\cos xt\;dx.$$

Тогда $F_c(t)\ge0,$  $F_c(t)\in L_1(0,\;\infty)$ и
$$
 \int\limits_0^\infty
F_c(t)dt=\sqrt\frac{\pi}{2}f(0). $$
\end{theorem}

{\sc Доказательство.} Так как
$$F_c(t)=-\sqrt\frac{2}{\pi}\frac{1}{t}\int\limits_0^\infty f'_+(x)\sin xt\;dx,$$
  то теорема~\ref{Skoryk-Theorem-S2} есть прямое следствие
  теоремы~\ref{Skoryk-Theorem-S1}, применяемой к функции $-f'_+(x).$

Теорема~\ref{Skoryk-Theorem-S2} совпадает с теоремой 124 из книги
Титчмарша \cite{Skoryk-Titchmarsh}. Мы считаем, что приводимое там
доказательство сложнее нашего.

Аналогично теореме~\ref{Skoryk-Theorem-S1} доказывается следующая
теорема.

\begin{theorem}\label{Skoryk-Theorem-S3} \hskip-2mm.
Пусть $\varphi(x)$ --- убывающая на полуоси $(0,\;\infty)$ функция
со сходящимся интегралом
$\displaystyle\int\limits_0^\infty\varphi(x)\;dx$ и пусть
$$g_m(t)=\int\limits_{\frac{\pi}{t}(m-\frac{1}{2})}^\infty \varphi(x)\cos xt\;dx=
\frac{1}{t}\int\limits_{\pi(m-\frac{1}{2})}^\infty\varphi\left(\frac{u}{t}\right)
\cos u\;du, \;\;m=1,2,\ldots$$

Тогда $(-1)^m g_m(t)$ есть положительная функция на полуоси
$(0,\;\infty)$  и выполняется равенство
$$\int\limits_0^\infty \frac{g_m(t)}{t}\;dt=\int\limits_{\pi(m-\frac{1}{2})}^\infty
\frac{\cos u}{u}\;du\int\limits_0^\infty\varphi(\tau)\;d\tau.$$
\end{theorem}

Докажем аналог теоремы~\ref{Skoryk-Theorem-S2} для синус
преобразования Фурье.

\begin{theorem}\label{Skoryk-Theorem-3alpha} \hskip-2mm.
Пусть $f(x)$ --- убывающая выпуклая бесконечно малая в
бесконечности функция на полуоси $[0,\;\infty).$ Тогда синус
преобразование Фурье функции $f(x)$
$$
F_s(t)=\sqrt\frac{2}{\pi}\int\limits_0^\infty f(x)\sin xt\;dx
\n{formula-06261}
$$
представляется в виде
$$F_s(t)=\sqrt\frac{2}{\pi}\;\frac{1}{t}f\left(\frac{\pi}{2t}\right)+\psi(t),$$
где $\psi\in L_1(0,\;\infty).$ Если дополнительно предположить,
что сходится интеграл $$\int\limits_0^\infty
\frac{f(x)-f(0)}{x}\;dx,$$ то справедлива формула
$$F_s(t)=\sqrt\frac{2}{\pi}\frac{f(0)}{t}+\psi_1(t),$$ где $\psi_1(t)\in L_1(0,\;\infty).$
\end{theorem}

{\sc Доказательство.} Из условий теоремы следует, что к интегралу
(\ref{formula-06261}) применима формула интегрирования по частям.
Поэтому
$$\displaystyle F_s(t)=-\sqrt\frac{2}{\pi}\;\frac{1}{t}f(x)\cos xt\Bigl|_0^\infty+\sqrt\frac{2}{\pi}
\;\frac{1}{t}\int\limits_0^\infty f'_+(x)\cos xt dx=$$
$$=\sqrt\frac{2}{\pi}\;\frac{f(0)}{t}+\sqrt\frac{2}{\pi}\;\frac{1}{t}\left(f\left(\frac{\pi}{2t}\right)
-f(0)\right)+\sqrt\frac{2}{\pi}\;\frac{1}{t}\int\limits_0^{\textstyle\frac{\pi}{2t}}
f'_+(x)(\cos xt-1)dx+$$
$$+\sqrt\frac{2}{\pi}\;\frac{1}{t}\int\limits_{\textstyle \frac{\pi}{2t}}^\infty f'_+(x)\cos xt dx.$$
Пусть $A_1(t)$ и $A_2(t)$ --- предпоследнее и последнее слагаемые
в написанной сумме. По теореме~\ref{Skoryk-Theorem-S3} имеем, что
$A_2\in L_1(0,\infty).$ Далее имеем
$$\int\limits_0^\infty |A_1(t)|dt\le\sqrt\frac{2}{\pi}\int\limits_0^\infty\frac{1}{t^2}
\int\limits_0^{\textstyle
\frac{\pi}{2}}\left|f'_+\left(\frac{u}{t}\right)\right|(1-\cos
u)du\;dt=$$
$$=\sqrt\frac{2}{\pi}\int\limits_0^{\textstyle\frac{\pi}{2}}(1-\cos u)\int\limits_0^\infty
\frac{1}{t^2}\left|f'_+\left(\frac{u}{t}\right)\right|dt\;du=\sqrt\frac{2}{\pi}
\int\limits_0^\infty\left|f'_+(\tau)\right|d\tau\int\limits_0^{\textstyle\frac{\pi}{2}}\frac{1-\cos
u}{u}\;du.$$

Из приведенных рассуждений следуют утверждения теоремы.

Дополнительно можно заметить, что
$$\int\limits_0^\infty\psi(t)dt=\sqrt\frac{2}{\pi}f(0)\left(\int\limits_0^{\textstyle\frac{\pi}{2}}
\frac{1-\cos u}{u}\;du
-\int\limits_{\textstyle\frac{\pi}{2}}^\infty\frac{\cos u}{u}\;
du\right).$$ Теорема доказана.

Заметим, что в отличие от теоремы~\ref{Skoryk-Theorem-S2},
теорема~\ref{Skoryk-Theorem-3alpha} отсутствует в книге Титчмарша.

Следующую теорему можно рассматривать как усиление той части
теоремы~\ref{Skoryk-Theorem-S2}, где говорится, что $F_c(t)\in
L_1(0,\;\infty).$ Мы заменяем требование выпуклости функции $f(x)$
на более слабое. Естественно, что при этом теряется свойство
положительности функции $F_c(t).$ Стандартный символ
$\displaystyle \bigvee\limits_a^b f$ обозначает вариацию функции
$f$ на промежутке $[a,\;b].$

\begin{theorem}\label{Skoryk-Theorem-d1} \hskip-2mm.
Пусть функция $R(x)$ удовлетворяет условиям:

1) $R(x)$ --- чётная функция,

2) $R(x)\to 0\;\;(x\to+\infty),$

3) $R(x)$ абсолютно непрерывна на любом сегменте $[0, a],\;\;a>0,$

4) $R'(x)\in L_1(0, \infty),$

5) существует число $b\ge0$ такое, что выполняются условия:

a) интеграл
$$
\int\limits_1^\infty \frac{|g(t)|}{t}\;dt<\infty,
\n{formula-07441}
$$
где $\displaystyle g(t)=\int\limits_0^b R'(x)\sin xt\;dx,$
сходится,

б) существует продолжение $R'$ на всю вещественную ось такое, что
функция $R'(x)$ имеет ограниченную вариацию на полуоси
$[b,\;\infty),$ причём
$\displaystyle\int\limits_b^\infty\bigvee\limits_x^\infty
R'\;dx<\infty.$

Тогда $R(x)\in\FR.$

\end{theorem}

{\sc Доказательство.}

Функция $R'(x)$ принадлежит $L_1(0,\;\infty)$ и имеет ограниченную
вариацию на полуоси $[b,\;\infty).$ Поэтому
$R'(x)\to0(x\to+\infty).$ При $x\ge b$ справедливо представление
$R'(x)=\varphi_1(x)-\varphi_2(x),$ где
$\varphi_1(x)=\bigvee\limits_x^\infty R',$
$\varphi_2(x)=\bigvee\limits_x^\infty R'-R'(x).$ Справедливы
соотношения $\varphi_1(x)\downarrow0,$ $\varphi_2(x)\downarrow0$
($x\to+\infty$), $\varphi_1,$ $\varphi_2\in L_1(b,\;\infty).$
Определяя $\varphi_1(x)=\varphi_1(b),$ $\varphi_2(x)=\varphi_2(b)$
при $x\in[0,\;b]$ мы получим, что функции $\varphi_1(x)$ и
$\varphi_2(x)$ убывают на полуоси $[0,\;\infty)$ и принадлежат
пространству $L_1(0,\;\infty).$ Обозначим
$$R_1(x)=-\int\limits_x^\infty \varphi_1(t)\;dt,\;\;
R_2(x)=-\int\limits_x^\infty
\varphi_2(t)\;dt,\;\;x\in[0,\;\infty)$$ и продолжим эти функции
чётным образом на всю ось $(-\infty,\;\infty).$ Обозначим также
$A(x)=R_1(x)-R_2(x).$ По теореме~\ref{Skoryk-Theorem-S2} функции
$R_1$ и $R_2$ принадлежат алгебре $\FR.$ Тогда $A\in\FR.$ Разность
$B(x)=R(x)-A(x)$ является непрерывной функцией на оси
$(-\infty,\;\infty),$ которая равна нулю при $|x|\ge b.$ На
сегменте $[0,\;b]$ функция $A(x)$ является линейной функцией.
Далее имеем
$$\int\limits_0^\infty B(t)\cos xt\;dt=-\frac{1}{t}\int\limits_0^\infty B'(t)\sin xt\;dt=
-\frac{1}{t}\int\limits_0^b R'(t)\sin
xt\;dt+\frac{a}{t}\int\limits_0^b \sin xt\;dt=$$
$$=-\frac{1}{t}\int\limits_0^bR'(t)\sin xt\;dt+a\frac{1-\cos
bt}{t^2}.$$ Теперь из условия 5) теоремы следует, что $\widehat
B\in L_1(1,\infty).$ Поскольку функция $\widehat B$ непрерывная и
чётная, то $\widehat B\in L_1(-\infty, \infty).$ Функция $B$
непрерывна и принадлежит $L_1(-\infty, \infty).$ Кроме того,
$\widehat B\in L_1(-\infty, \infty).$ Поэтому $B\in\FR,$
$R\in\FR.$

Теорема доказана.

Заметим, что неравенство (\ref{formula-07441}) можно заменить на
более сильное неравенство $\displaystyle\int\limits_0^b
\left|R'(x)\right|^p\;dx<\infty$ с $p>1.$ Можно считать, что
$p\in(1,\;2].$ В этом случае $g(t)\in L_q(0,\;\infty),$
$\displaystyle \frac{1}{p}+\frac{1}{q}=1.$ И тогда неравенство
(\ref{formula-07441}) следует из неравенства Шварца.

Далее мы сформулируем аналог теоремы~\ref{Skoryk-Theorem-d1} для
нечётных функций.

\begin{theorem}\label{Skoryk-Theorem-d2}\hskip-2mm.
Пусть функция $R(x)$ удовлетворяет условиям:

1) $R(x)$ --- нечётная функция, $R(0)=0,$

2) $R(x)\to 0\;\;(x\to+\infty),$

3) интеграл
$\displaystyle\int\limits_0^\infty\frac{|R(x)|}{x}\;dx$ сходится,

4) $R(x)$ абсолютно непрерывна на любом сегменте $[0, a],\;\;a>0,$

5) $R'(x)\in L_1(0, \infty),$

6) существует число $b\ge0$ такое, что выполняются условия:

a) $\displaystyle \int\limits_1^\infty
\frac{|g(t)|}{t}\;dt<\infty,$ где $\displaystyle
g(t)=\int\limits_0^b R'(x)\cos xt\;dx,$

б) существует продолжение $R'$ на всю вещественную ось такое, что
функция $R'(x)$ имеет ограниченную вариацию на полуоси
$[b,\;\infty),$ причём
$$\displaystyle\int\limits_b^\infty\bigvee\limits_x^\infty
R'\;dx<\infty.$$

Тогда $R(x)\in\FR.$

\end{theorem}

Доказательство теоремы~\ref{Skoryk-Theorem-d2} проводится по той
же схеме, что и доказательство теоремы~\ref{Skoryk-Theorem-d1},
только ссылку на теорему~\ref{Skoryk-Theorem-S2} нужно заменить
ссылкой на теорему~\ref{Skoryk-Theorem-3alpha}.

В однотипных теоремах~\ref{Skoryk-Theorem-d1} и
\ref{Skoryk-Theorem-d2}  рассматривались случаи, когда $R(x)$
является чётной функцией (теорема~\ref{Skoryk-Theorem-d1}) и
нечётной функцией (теорема~\ref{Skoryk-Theorem-d2}). Из этих
теорем легко получается аналогичная теорема, где отсутствуют
ограничения на чётность функции $R(x).$ Эта теорема доказывается
применением теорем~\ref{Skoryk-Theorem-d1} и
\ref{Skoryk-Theorem-d2} к функциям
$\displaystyle\frac{1}{2}\left(R(x)+R(-x)\right)$ и
$\displaystyle\frac{1}{2}\left(R(x)-R(-x)\right).$

\begin{theorem}\label{Skoryk-Theorem-6alpha}\hskip-2mm.
Пусть функция $R(x)$ удовлетворяет условиям:

1) $R(x)\to0\;\;(x\to\pm\infty),$

2) интеграл
$\displaystyle\int\limits_0^\infty\frac{|R(x)-R(-x)|}{x}\;dx$
сходится,

3) $R(x)$ абсолютно непрерывна на любом сегменте,

4) $R'(x)\in L_1(-\infty,\;\infty),$

5) существует число $b\ge0$ такое, что выполняются условия:

a)
$$
\int\limits_1^\infty\frac{|g(t)|}{t}\;dt<\infty,\quad\int\limits_1^\infty\frac{|h(t)|}{t}\;dt<\infty,\n{formula-07442}
$$
где
$$g(t)=\int\limits_0^b\left(R'(x)-R'(-x)\right)\sin xt \;dx,\quad
h(t)=\int\limits_0^b\left(R'(x)+R'(-x)\right)\cos xt \;dx,$$

б) существует продолжение $R'$ на всю вещественную ось такое, что
выполняются условия
$$\int\limits_{-\infty}^{-b}\bigvee_{-\infty}^xR'\;dx<\infty,\qquad\int\limits_b^\infty\bigvee_x^\infty R'\;dx<\infty.$$

Тогда $R(x)\in\FR.$
\end{theorem}

Справедливо замечание, аналогичное замечанию к
теореме~\ref{Skoryk-Theorem-d1}. Неравенство $$\int\limits_{-b}^b
\left|R'(x)\right|^p\;dx<\infty$$ с некоторым $p>1$ сильнее
неравенств (\ref{formula-07442}).

В следующей теореме рассматривается случай, когда функция $R'(x)$
кусочно монотонна и интегрируемость функции $|R'(x)|^p$ может
нарушаться в конечном числе точек.

\begin{theorem}\label{Skoryk-Theorem-C0} \hskip-2mm.
Пусть $R(x)$ --- непрерывная стремящаяся к нулю при
$x\to\pm\infty$ функция и пусть существует такое покрытие
$$(-\infty,\;a_1],\;\;[a_1,\;a_2],\ldots,\;\;[a_{n-1},\;a_n],\;\;[a_n,\;\infty)$$
вещественной оси, что на каждом из выписанных выше множеств одна
из функций $R(x)$ и $-R(x)$ является выпуклой.

Тогда для того, чтобы $R(x)\in\FR,$ необходимо и достаточно, чтобы
сходились интегралы
$$\int\limits_1^\infty\frac{|R(x)-R(-x)|}{x}\;dx,\quad\int\limits_0^1\frac{|R(a_k+x)-R(a_k-x)|}{x}\;dx,
 \quad k=1,\ldots,n.$$

\end{theorem}

{\sc Доказательство.}
 Вначале построим специальную систему функций \linebreak $\bigl\{h_k(x)\bigr\}_{k=0}^n.$

 Выберем число $c>0$ таким, чтобы выполнялись условия:
 $c\in(a_n+5,\infty),$ $-c\in(-\infty,a_1-5),$ точки $c$ и
 $-c$ --- точки дифференцируемости функции $R(x).$ Определим
 функцию $H_0(x)$  на полуосях $(-\infty,\;-c],$ $[c,\;\infty)$
 так, чтобы выполнялись условия:

 1) на каждой из этих полуосей функция $H_0(x)$ имеет то же направление
 выпуклости, что и $R(x),$

 2) $H_0(-c)=R(-c),$ $H'_0(-c)=R'(-c),$ $H_0(c)=R(c),$ $H'_0(c)=R'(c),$

 3) функция $H_0(x)$ является дважды непрерывно дифференцируемой
 на каждой из полуосей $(-\infty,-c]$ и $[c,\infty),$

 4) функции $H_0(x),$ $H'_0(x),$ $H''_0(x)$ принадлежат пространствам
 $L_1(-\infty, -c)$ и $L_1(c, \infty).$

 Функцию $h_0(x)$ определим формулой
 $$h_0(x)=\left\{\begin{array}{cl}
\displaystyle R(x)-H_0(x),&\quad |x|\ge c,
\\[5pt]
0,&\quad |x|<c.
\end{array}\right.$$

Далее выберем число $\delta>0$ так, чтобы выполнялись условия :

1) $a_k+2\delta<a_{k+1}-2\delta,\quad k=1,\ldots,n{-}1,$

2) точки $a_k\pm\delta$  были точками дифференцируемости функции
$R(x),$

3) на каждом из интервалов $(a_k-\delta, a_k),$ $(a_k,
a_k+\delta),$ $k=1,\ldots,n$ функция $R'_+(x)$ не меняет знак.

Далее для $k=1,\ldots,n$ функции $h_k(x)$ определяем следующим
образом
$$h_k(x)=\left\{\begin{array}{cl}
0,& x\in(-\infty, a_k-2\delta], \\[3pt]
P_k(x),& x\in[a_k-2\delta, a_k-\delta],\\[3pt]
R(x),& x\in[a_k-\delta, a_k+\delta],\\[3pt]
Q_k(x),& x\in[a_k+\delta, a_k+2\delta],\\[3pt]
0,& x\in[a_k+2\delta,\infty),
\end{array}\right.$$
где $P_k(x)$ и $Q_k(x)$ --- многочлены не выше третьей степени,
выбираемые так, чтобы функция $h_k(x)$ была дифференцируемой в
точках $x=a_k\pm\delta,$ $x=a_k\pm2\delta.$

Обозначим
$$h(x)=\sum\limits_{k=0}^n h_k(x),\quad \psi(x)=R(x)-h(x).$$
Функция $\psi(x)$ обладает свойствами:

1) $\psi(x)$ --- непрерывная функция на всей оси и $\psi(x)\to0$
при $x\to\pm\infty,$

2) $\psi(x)\in L_1(-\infty,\;\infty),$

3) у функции $\psi(x)$ всюду существует правая производная
$\psi'_+(x)$ и эта производная является функцией ограниченной
вариации на оси $(-\infty,\;\infty).$

Из этих свойств легко следует, что преобразование Фурье
$\widehat\psi(t)$ функции $\psi(x)$ принадлежит пространству
$L_1(-\infty,\infty),$ а сама функция $\psi(x)$ принадлежит
алгебре $\FR.$

Далее мы будем исследовать свойства преобразований Фурье функций
$h_k(x).$ Пусть $k=1,\ldots,n.$ Имеем
$$\widehat h_k(t)=\frac{1}{\sqrt{2\pi}}\int\limits_{-\infty}^\infty h_k(x)e^{-ixt}\;dx=
\frac{e^{-ia_kt}}{\sqrt{2\pi}}\int\limits_{-\infty}^\infty
h_k(a_k+u)e^{-iut}\;du=$$
$$=\frac{e^{-ia_kt}}{\sqrt{2\pi}}\left(\int\limits_0^\infty\left(h_k(a_k+u)+
h_k(a_k-u)\right)\cos ut\;du-\right.$$
$$\left.-i\int\limits_0^\infty\left(h_k(a_k+u)-h_k(a_k-u)\right)\sin
ut\;du\right).$$

Обозначим
$$c_k(t)=\int\limits_0^\infty\bigl(h_k(a_k+u)+h_k(a_k-u)\bigr)\cos ut\;du,$$
$$d_k(t)=\int\limits_0^\infty\bigl(h_k(a_k+u)-h_k(a_k-u)\bigr)\sin ut\;du.$$
Далее находим
$$c_k(t)=-\frac{1}{t}\int\limits_0^\infty\bigl(h'_k(a_k+u)-h'_k(a_k-u)\bigr)\sin
ut\;du=$$ $$=-\frac{1}{t}\int\limits_0^\delta h'_k(a_k+u)\sin
ut\;du+\frac{1}{t}\int\limits_0^\delta h'_k(a_k-u)\sin ut\;du-$$
$$-\frac{1}{t}\int\limits_\delta^\infty
\bigl(h'_k(a_k+u)-h'_k(a_k-u)\bigr)\sin
ut\;du=c_{k1}(t)+c_{k2}(t)+c_{k3}(t).$$

Функция $h_k$  в общем случае не дифференцируема. В тех точках,
где производная не существует, под производной следует понимать
правую производную. Поскольку функция $h'_k(a_k+u)-h'_k(a_k-u)$
имеет ограниченную вариацию на полуоси $[\delta,\;\infty),$ то
интегрирование по частям даёт существование  постоянной $M_{k3}$
такой, что выполняется неравенство
$$|c_{k3}(t)|\le\frac{M_{k3}}{t^2}.$$

Рассмотрим функцию $c_{k1}(t).$ На интервале $(0,\delta)$
выполняется равенство $h'_k(a_k+u)=R'(a_k+u).$ Если $R'(a_k+0)$
--- конечная величина, то вновь интегрирование по частям приведёт
к неравенству
$$|c_{k1}(t)|\le\frac{M_{k1}}{t^2}.$$

Пусть теперь $|R'(a_k+0)|=\infty.$ Не ограничивая общности, можно
считать, что $R'(a_k+0)=+\infty$ (иначе нужно рассматривать
функцию $-c_{k1}(t)$). В силу ограничений, наложенных на $\delta,$
имеем $R'(a_k+\delta)\ge0.$ Поэтому если функцию $h'_k(a+u)$
продолжить с полуинтервала $(0,\delta]$ на полуось $[0,\infty),$
полагая продолженную функцию равной $0$ на полуоси $(\delta,
\infty),$ то продолженная функция будет убывающей на полуоси $(0,
\infty).$ В рассматриваемом случае применение
теоремы~\ref{Skoryk-Theorem-S1} даёт соотношение $c_{k1}(t)\in
L_1(0, \infty).$ Однако, в любом случае $c_{k1}(t)\in L_1(1,
\infty).$ Аналогично можно доказать, что $c_{k2}(t)\in L_1(1,
\infty).$ Из сказанного следует, что $c_k(t)\in L_1(1, \infty).$

Переходим к исследованию $d_k(t).$ Имеем при $\displaystyle
t>\frac{\pi}{2\delta}$
$$d_k(t)=\frac{1}{t}\int\limits_0^\infty\bigl(h'(a_k{+}u)+h'(a_k{-}u)\bigr)\cos ut\;du=
\frac{1}{t}\left(h_k\bigl(a_k{+}\frac{\pi}{2t}\bigr){-}h_k\bigl(a_k{-}\frac{\pi}{2t}\bigr)\right)-$$
$$-\frac{1}{t}\int\limits_0^\frac{\pi}{2t}\left(h'(a_k+u)+h'(a_k-u)\right)(1-\cos ut)\;du
+\frac{1}{t}\int\limits_{\textstyle \frac{\pi}{2t}}^\delta
h'_k(a_k+u)\cos ut\;du+$$
$$+\frac{1}{t}\int\limits_{\textstyle \frac{\pi}{2t}}^\delta h'_k(a_k-u)\cos ut\;du+\frac{1}{t}
\int\limits_\delta^\infty\left(h'_k(a_k+u)+h'_k(a_k-u)\right)\cos
ut\;du=$$
$$=\frac{1}{t}\left(h_k\bigl(a_k+\frac{\pi}{2t}\bigr)-h_k\bigl(a_k-\frac{\pi}{2t}\bigr)\right)+d_{k1}(t)+
d_{k2}(t)+d_{k3}(t)+d_{k4}(t).$$

Так как функция $h'_k(a_k{+}u)+h'_k(a_k{-}u)$ имеет ограниченную
вариацию на полуоси $[\delta,\infty),$ то интегрирование по частям
приводит к оценке
$$|d_{k4}(t)|\le\frac{M_{k4}}{t^2}$$
с некоторой постоянной $M_{k4}.$

Если обозначить $\varphi(u)=|h'_k(a_k+u)+h'_k(a_k-u)|,$  то
получим
$$\int\limits_0^\infty|d_{k1}(t)|dt\le\int\limits_0^\infty\frac{1}{t}\int\limits_0^
{\textstyle \frac{\pi}{2t}}\varphi(u)(1-\cos ut)du\;dt=
\int\limits_0^\infty\frac{1}{t^2}\int\limits_0^{\textstyle
\frac{\pi}{2}}\varphi\left(\frac{v}{t}\right) (1-\cos v)dv\;dt=$$
$$=\int\limits_0^{\textstyle \frac{\pi}{2}}(1-\cos
v)\int\limits_0^\infty\frac{1}{t^2}\varphi\left(\frac{v}{t}\right)dt\;dv=\int\limits_0^\infty\varphi(\tau)d\tau\int\limits_0^{\textstyle
\frac{\pi}{2}}\frac {1-\cos v}{v}\;dv.$$

Оценка функции $d_{k2}(t)$ проводится аналогично оценке функции
$c_{k1}(t),$ только ссылку на теорему~\ref{Skoryk-Theorem-S1}
нужно заменить ссылкой на теорему~\ref{Skoryk-Theorem-S3}. Это
даёт $d_{k2}(t)\in L_1(\frac{\pi}{\delta},\;\infty).$ Аналогично
$d_{k3}(t)\in L_1(\frac{\pi}{\delta},\;\infty).$

Таким образом, при $\displaystyle t\ge\frac{\pi}{\delta}$
выполняется равенство
$$d_k(t)=\frac{1}{t}\left(R\bigl(a_k+\frac{\pi}{2t}\bigr)-R\bigl(a_k-\frac{\pi}{2t}\bigr)\right)+
\widetilde{d_k}(t),$$  где $\tilde{d_k}(t)\in
L_1(\frac{\pi}{\delta},\;\infty).$ Возвращаясь к функции $\widehat
h_k(t)$ получим, что при $\displaystyle t\ge\frac{\pi}{\delta}$
выполняется равенство
$$\widehat h_k(t)=\frac{e^{-ia_kt}}{\sqrt{2\pi}i}\frac{1}{t}\left(R\bigl(a_k+\frac{\pi}{2t}\bigr)-
R\bigl(a_k-\frac{\pi}{2t}\bigr)\right)+g_k(t),$$ где $g_k(t)\in
L_1(\frac{\pi}{\delta},\;\infty).$

Осталось исследовать функцию
$$\widehat h_0(t)=\frac{1}{\sqrt{2\pi}}\int\limits_{-\infty}^\infty h_0(x)e^{-ixt}\;dx=
\frac{1}{\sqrt{2\pi}ti}\int\limits_{-\infty}^\infty
h'_0(x)e^{-ixt}\;dx=$$
$$=\frac{1}{\sqrt{2\pi}ti}\left(\int\limits_0^\infty\bigl(h'_0(x)+h'_0(-x)\bigr)
\cos xt\;dx-i\int\limits_0^\infty\bigl(h'_0(x)-h'_0(-x)\bigr)\sin
xt\;dx\right).$$ Обозначим
$$c_0(t){=}\int\limits_0^\infty\bigl(h'_0(x){+}h'_0(-x)\bigr)\cos
xt\;dx{=}\int\limits_c^\infty\bigl(R'(x)+R'(-x)-H'_0(x)-H'_0(-x)\bigr)\cos
xt\;dx,$$
$$d_0(t){=}\int\limits_0^\infty\bigl(h'_0(x)-h'_0(-x)\bigr)\sin
xt\;dx{=}\int\limits_c^\infty\bigl(R'(x)-R'(-x)-H'_0(x)+H'_0(-x)\bigr)\sin
xt\;dx.$$

Пусть $f(x)$ --- функция, определённая на полуоси $[c, \infty).$
Обозначим через $\widetilde{f}(x)$ продолжение этой функции на
полуось $[0, \infty),$ причём $\widetilde{f}(x)=f(c)$ при $x\in[0,
c].$ Справедливы равенства
$$c_0(t)=\int\limits_0^\infty\left(\widetilde R'(x)+\widetilde R' (-x)-\widetilde H'_0 (x)-\widetilde H'_0
(-x)\right)\cos xt\;dx,$$
$$d_0(t)=\int\limits_0^\infty\left(\widetilde R'(x)-\widetilde R'(-x)-\widetilde H'_0(x)+\widetilde H'_0
(-x)\right)\sin xt\;dx=$$
$$=d_{01}(t)+d_{02}(t)+d_{03}(t)+d_{04}(t).$$

Функция $\widetilde R'(x)$ монотонна на полуоси $[0,\infty)$ и
принадлежит пространству $L_1(0, \infty).$ Поэтому из
теоремы~\ref{Skoryk-Theorem-S1} следует сходимость интеграла
$\displaystyle\int\limits_0^\infty\frac{|d_{01}(t)|}{t}\;dt.$ Эти
же рассуждения применимы для функций $d_{02},$ $d_{03},$ $d_{04}.$
Тем самым доказано, что
$$\int\limits_0^\infty\frac{|d_0(t)|}{t}\;dt<\infty.$$

Обозначим $\varphi(x)=\widetilde R'(x)+\widetilde
R'(-x)-\widetilde H'_0(x)-\widetilde H'_0(-x).$ Имеем

$$c_0(t)=\int\limits_0^{\textstyle \frac{\pi}{2t}}\varphi(x)dx-\int\limits_0^{\textstyle \frac{\pi}{2t}}
\varphi(x)(1-\cos xt)dx+ \int\limits_{\textstyle
\frac{\pi}{2t}}^\infty\widetilde R'(x)\cos xt\;dx+$$
$$+\int\limits_{\textstyle \frac{\pi}{2t}}^\infty\widetilde R'(-x)\cos xt\;dx-
\int\limits_{\textstyle \frac{\pi}{2t}}^\infty\widetilde
H'_0(x)\cos xt\;dx-\int\limits_{\textstyle
\frac{\pi}{2t}}^\infty\widetilde H'_0(-x)\cos xt\;dx
=\sum\limits_{k=1}^6 c_{0k}(t).$$

Справедливо неравенство
$$\int\limits_0^\infty\frac{|c_{02}(t)|}{t}\;dt\le\int\limits_0^\infty\frac{1}{t}
\int\limits_0^{\textstyle \frac{\pi}{2t}}|\varphi(x)|(1-\cos
xt)dx\;dt=$$
$$=\int\limits_0^\infty\frac{1}{t^2}\int\limits_0^{\textstyle \frac{\pi}{2}}|\varphi(\frac{u}{t})|(1-\cos u)du\;dt
=\int\limits_0^\infty|\varphi(v)|dv\int\limits_0^{\textstyle
\frac{\pi}{2}}\frac{1-\cos u}{u}\;du.$$

Функция $\widetilde R'(x)$ монотонна на полуоси $(0,\infty)$ и
принадлежит пространству $L_1(0, \infty).$ По
теореме~\ref{Skoryk-Theorem-S3}
$$\int\limits_0^\infty\frac{|c_{03}(t)|}{t}\;dt<\infty.$$

Аналогичные рассуждения справедливы для функций $c_{04}(t),$
$c_{05}(t),$ $c_{06}(t).$

Если $\displaystyle t\ge\frac{\pi}{2c},$ то $c_{01}(t)=0.$ В
противном случае
$$c_{01}(t)=\int\limits_c^{\textstyle \frac{\pi}{2t}}\left(\widetilde R'(x)+\widetilde R'(-x)
-\widetilde H'_0(x)-\widetilde H'_0(-x)\right)dx=$$
$$=R\left(\frac{\pi}{2t}\right)-R\left(-\frac{\pi}{2t}\right)-H_0\left(\frac{\pi}{2t}\right)+H_0\left(-\frac{\pi}{2t}\right)$$

Если обозначить

$$R_1(t)=\left\{\begin{array}{cl}
\displaystyle
R\left(\frac{\pi}{2t}\right)-R\left(-\frac{\pi}{2t}\right),&\displaystyle
t\in\left(0,\;\frac{\pi}{2c}\right),
\\[10pt]
0,& \displaystyle t\ge\frac{\pi}{2c},
\end{array}\right.$$
$$H_1(t)=\left\{\begin{array}{cl}
\displaystyle
H_0\left(\frac{\pi}{2t}\right)-H_0\left(-\frac{\pi}{2t}\right),&
\displaystyle t\in\left(0,\;\frac{\pi}{2c}\right),
\\[10pt]
0,& \displaystyle  t\ge\frac{\pi}{2c},
\end{array}\right.$$
то будет выполняться равенство $c_{01}(t)=R_1(t)-H_1(t).$ Легко
проверяется неравенство
$$\int\limits_0^\infty\frac{|H_1(t)|}{t}\;dt<\infty.$$

Из сказанного следует, что справедливо представление
$$\widehat h_0(t)=\frac{1}{\sqrt{2\pi}\;it}R_1(t)+g_0(t),$$
где $g_0(t)\in L_1(0, \infty).$

Теперь легко заканчивается доказательство теоремы. Если написанные
в условии теоремы интегралы сходятся, то функции $\widehat
h_0(t),$ $\widehat h_1(t),$ $\ldots,$ $\widehat h_n(t)$
принадлежат $L_1(0, \infty),$ а значит и $L_1(-\infty,\infty).$
Нужно иметь в виду, что при $k=1,\ldots,n$ функции $\widehat
h_k(t)$ непрерывны и поэтому принадлежат
$L_1(0,\;\frac{\pi}{\delta}).$ Из этого следует, что $\widehat
R(t)\in L_1(-\infty,\infty).$ Имеем
$$\widehat R(t)=\frac{1}{\sqrt{2\pi}}\int\limits_{-\infty}^\infty
R(x)e^{-itx}\;dx.$$ Если написанный интеграл рассматривать как
несобственный интеграл с особыми точками $\pm\infty,$ то из
формулы интегрирования по частям следует, что функция $\widehat
R(t)$ корректно определена для всех $t$ за возможным исключением
$t=0.$ Кроме того, $R(x)\to0$ при $x\to\pm\infty.$ Локальной
интегрируемости функции $\widehat R(t)$ достаточно для
справедливости формулы обращения
$$R(x)=\frac{1}{\sqrt{2\pi}}\lim\limits_{\lambda\to\infty}\int\limits_{-\lambda}^\lambda
\left(1-\frac{|t|}{\lambda}\right)\widehat R(t)e^{itx}\;dt$$ почти
для всех $x$ (это следует из теоремы 113
\cite{Skoryk-Titchmarsh}). В нашем случае, поскольку $\widehat
R(t)\in L_1(-\infty,\infty),$ то это равенство можно переписать в
виде
$$R(x)=\frac{1}{\sqrt{2\pi}}\int\limits_{-\infty}^\infty \widehat R(t)e^{itx}\;dt.$$
Так как обе части равенства являются непрерывными функциями, то
оно выполняется для всех вещественных $x.$ Тем самым $R\in\FR.$

Пусть теперь $R(x)\in\FR.$ Тогда $\displaystyle
h(x)=\sum\limits_{k=0}^n h_k(x)\in\FR.$ Так как справедливо
неравенство  $\rho({\rm supp}\,h_k,\,{\rm supp}\,h_j)>0$ при $k\ne
j,$ то каждая из функций $h_k(x)\in\FR.$ Поэтому
$$h_k(x)=\frac{1}{\sqrt{2\pi}}\int\limits_{-\infty}^\infty \varphi_k(t)e^{-itx}\;dt,$$
где $\varphi_k\in L_1(-\infty, \infty).$ По формуле обращения
$$\varphi_k(t)=\lim\limits_{\varepsilon\to0}\frac{1}{\sqrt{2\pi}}\int\limits
_{-\infty}^\infty h_k(x)e^{-\varepsilon
x^2}e^{ixt}\;dx=\frac{1}{\sqrt{2\pi}}\int\limits _{-\infty}^\infty
h_k(x)e^{ixt}dx.$$ Последнее равенство для $k=1,\ldots,n$
обосновывается с помощью теоремы Лебега о предельном переходе. В
случае $k=0$ оно следует из равномерной сходимости
соответствующего несобственного интеграла на множестве
$\{\varepsilon: 0<\varepsilon<\infty\}$ (применяется признак Абеля
равномерной сходимости несобственных интегралов).

Таким образом, $\widehat h_k(t)\in L_1(0, \infty).$ Из этого в
свою очередь следует сходимость интегралов выписанных в условии
теоремы.

Теорема доказана.

В качестве примера рассмотрим функцию
$$R(x)=\frac{1}{\sqrt{1+\ln^2 |x|}}\;.$$
Имеем
$$R'(x)=\frac{-\ln|x|}{x\left(1+\ln^2|x|\right)^{\textstyle\frac{3}{2}}}, \quad R''(x)=
\frac{\ln^3|x|+2\ln^2|x|+\ln|x|-1}{x^2\left(1+\ln^2|x|\right)^{\textstyle\frac{5}{2}}}\;.$$
Для этой функции выполняются все условия
теоремы~\ref{Skoryk-Theorem-C0}. Поэтому $R(x)\in\FR.$ Отметим
ещё, что для функции
$$R(x)=\frac{1}{\ln_k(a_k+|x|)},$$ где $a_k\ge
e_k,$ ($e_1=e,$ $e_{k+1}=e^{e_k}$) проверка условий
теоремы~\ref{Skoryk-Theorem-C0} не вызывает никаких затруднений.
Поэтому $R(x)\in\FR.$

 \sect{Свойства функций из алгебры  $\KR$}\label{razdel_KR}

Известны сложности, связанные с описанием функций из алгебры
$\FR,$ а тем более из алгебры $\KR$. Об этом, в частности,
написано в работе \cite[глава 1, раздел~6]{Skoryk-Gurariy},
Конечно, в пунктах 12 и 13 из вступления приведены критерии
принадлежности функций алгебрам $\FR$ и $\KR.$ Однако, авторам
неизвестны случаи эффективного применения этих критериев. В
современных руководствах по гармоническому анализу эти критерии
зачастую не приводятся.

В предыдущем параграфе приведены достаточные условия
принадлежности функции алгебре $\FR.$ Далее приводятся некоторые
новые свойства функций из этих алгебр. Тем самым даются
необходимые условия вхождения функций в указанные алгебры.

\begin{theorem}\label{Skoryk-Theorem-3} \hskip-2mm. Пусть функция
$\varphi\in\KR$ и нечётная. Тогда для любого невещественного $z$
несобственный интеграл
$$
\int\limits_0^\infty\frac{\varphi(\lambda)}{\lambda+z}\,d\lambda
\n{Skoryk-formula-3.1}
$$
сходится.
\end{theorem}

{\sc Доказательство.} Из условий теоремы следует, что
$$\varphi(\lambda)=\int\limits_0^\infty \sin \lambda t \,d\mu(t),$$
где $\mu\in M({\Bbb R}).$ Поскольку величина
$$\int\limits_0^N \left(\int\limits_0^\infty \left|\frac{\sin \lambda t}{\lambda+z}\right|
\,d|\mu|(t) \right)d\lambda $$ является конечной, то по теореме
Тоннели \cite[гл.~III, \S~11, пункт~4]{Skoryk-Dandorf} функция
$\displaystyle\frac{\sin \lambda t}{\lambda+z}$ принадлежит
пространству $L_1$ на множестве $[0, N]{\times}[0, \infty)$ по
мере $d\lambda\times d\mu(t).$ Из теоремы Фубини следует, что
$$I_N(z)=\int\limits_0^N \left(\int\limits_0^\infty \frac{\sin
\lambda t}{\lambda+z} \,d\mu(t)
\right)d\lambda=\int\limits_0^\infty \left(\int\limits_0^N
\frac{\sin \lambda t}{\lambda+z} \,d\lambda \right)d\mu(t).$$
Имеем
$$A(N,t)=\int\limits_0^N \frac{\sin \lambda t}{\lambda+z}
\,d\lambda =\int\limits_0^N \frac{\sin \lambda t}{\lambda}
\,d\lambda-z\int\limits_0^N \frac{\sin \lambda
t}{\lambda(\lambda+z)} \,d\lambda=A_1(N,t)+A_2(N,t).$$ Из
равенства\vskip-4mm
$$A_1(N,t)=\int\limits_0^{Nt} \frac{\sin \lambda }{\lambda}
\,d\lambda$$ следует ограниченность функции $A_1(N,t)$ на
множестве $[0,\infty){\times}[0,\infty).$

Рассмотрим функцию
$$\frac{1}{\lambda+z}=\frac{\lambda+x-iy}{(\lambda+x)^2+y^2}.$$
Функция $\displaystyle \frac{\lambda+x}{(\lambda+x)^2+y^2},$ как
функция переменной $\lambda$ на оси $(-\infty,\infty),$ меняет
направление монотонности в точках $-x-y$ и $-x+y,$ а функция
$\displaystyle \frac{y}{(\lambda+x)^2+y^2}$ меняет направление
монотонности в точке $\lambda=-x.$

Если на сегменте $[\lambda_1,\lambda_2]$ функция $\displaystyle
\frac{\lambda+x}{(\lambda+x)^2+y^2}$ является монотонной, то по
второй теореме о среднем значении имеем
$$\int\limits_{\lambda_1}^{\lambda_2}\frac{(\lambda{+}x)\sin\lambda
t}{\lambda\left((\lambda{+}x)^2{+}y^2\right)}\,d\lambda=\frac{\lambda_1{+}x}{(\lambda_1{+}x)^2{+}y^2}
\int\limits_{\lambda_1}^\xi \frac{\sin\lambda
t}{\lambda}\,d\lambda{+}
\frac{\lambda_2{+}x}{(\lambda_2{+}x)^2{+}y^2}
\int\limits_{\xi}^{\lambda_2} \frac{\sin\lambda
t}{\lambda}\,d\lambda.$$ Для любого вещественного $\lambda$
выполняется неравенство $$\left|
\frac{\lambda+x}{(\lambda+x)^2+y^2}\right|\le \frac{1}{2|y|}\;.$$
Существует $M>0$ такое, что для любых $a,$ $b$ и $t$ выполняется
неравенство
$$\left|\int\limits_a^b \frac{\sin\lambda
t}{\lambda}\,d\lambda\right|\le M.$$ Заметим еще, что $$\left|
\frac{y}{(\lambda+x)^2+y^2}\right|\le \frac{1}{|y|}.$$ Теперь
нетрудно увидеть, что для функции $A(N,t)$ справедлива оценка
$$|A(N,t)|\le M\left(1+\frac{|z|}{|y|}\right).$$
Несобственный интеграл
$\displaystyle\int\limits_0^\infty\frac{\sin\lambda
t}{\lambda+z}\,d\lambda$ является сходящимся для любых $t\ge 0.$

Теперь из теоремы о мажорируемой сходимости следует, что функция
$I_N(z)$ имеет предел при $N\to\infty$ и этот предел равен
$$\int\limits_0^\infty \left(\int\limits_0^\infty \frac{\sin \lambda t}{\lambda+z}
\,d\lambda \right)d\mu(t). $$

Тем самым установлены не только сходимость интеграла
(\ref{Skoryk-formula-3.1}), но и выполнение равенства
$$\int\limits_0^\infty\frac{\varphi(\lambda)}{\lambda+z}\,d\lambda=
\int\limits_0^\infty \left(\int\limits_0^\infty \frac{\sin \lambda
t}{\lambda+z} \,d\lambda \right)d\mu(t). $$

Теорема доказана.

\begin{remark}\label{r1} \hskip-2mm.
B \cite[гл.~1, \S~4]{Skoryk-Stein} доказано более слабое
утверждение об ограниченности по переменной $b$ интеграла
$\displaystyle\int\limits_1^b\frac{\varphi(\lambda)}{\lambda}\,d\lambda$
для нечётных функций $\varphi$ из алгебры $\FR.$
\end{remark}

Заметим, что и в формулировке теоремы~\ref{Skoryk-Theorem-3}, и в
её доказательстве промежуток интегрирования $(0, \infty)$ можно
заменить на $(-\infty, 0).$ Из этого следует такая теорема.

\begin{theorem}\label{Skoryk-Theorem-4} \hskip-2mm. Пусть функция
$\varphi\in\KR$ и нечётная, а $z$ --- невещественное число. Тогда
сходится  несобственный интеграл $$
\int\limits_{-\infty}^\infty\frac{\varphi(\lambda)}{\lambda+z}\,d\lambda.
$$
Более того, если
$$\varphi(\lambda)=\int\limits_0^\infty \sin \lambda t\,
d\mu(t),\quad \mu\in M({\Bbb R}),$$ то выполняется равенство
$$\int\limits_{-\infty}^\infty\frac{\varphi(\lambda)}{\lambda+z}\,d\lambda=
\int\limits_0^\infty \left(\int\limits_{-\infty}^\infty \frac{\sin
\lambda t}{\lambda+z} \,d\lambda \right)d\mu(t). $$
\end{theorem}

\begin{theorem}\label{Skoryk-Theorem-5} \hskip-2mm. Если функция
$\varphi\in\KR$ и чётная, а $z$ --- невещественное число, то
существует
$$V.P.\int\limits_{-\infty}^\infty\frac{\varphi(\lambda)}{\lambda+z}\,d\lambda.
$$
Более того, если
$$
\varphi(\lambda)=\int\limits_0^\infty \cos \lambda t\,
d\mu(t),\quad \mu\in M(\Bbb R),\n{Skoryk-formula-4}$$ то
выполняется равенство
$$V.P. \int\limits_{-\infty}^\infty\frac{\varphi(\lambda)}{\lambda+z}\,d\lambda=
\int\limits_0^\infty \left(V.P.\int\limits_{-\infty}^\infty
\frac{\cos \lambda t}{\lambda+z} \,d\lambda \right)d\mu(t). $$
\end{theorem}

{\sc Доказательство.} Так как функция $\varphi\in\KR$ и чётная, то
равенство (\ref{Skoryk-formula-4}) выполняется для некоторой меры
$\mu\in M(\Bbb R).$

Рассмотрим
$$I_N(z)=\int\limits_{-N}^N
\left(\int\limits_0^\infty \frac{\cos \lambda t}{\lambda+z}
\,d\mu(t) \right)d\lambda.$$ Повторение соответствующих
рассуждений из доказательства теоремы~\ref{Skoryk-Theorem-3} даёт
$$I_N(z)=\int\limits_{0}^\infty
\left(\int\limits_{-N}^N \frac{\cos \lambda t}{\lambda+z}
\,d\lambda \right)d\mu(t).$$ Так как
$$A(N,t)=\int\limits_{-N}^N \frac{\cos \lambda t}{\lambda+z}
\,d\lambda =2z\int\limits_0^N \frac{\cos \lambda
t}{(z-\lambda)(z+\lambda)} \,d\lambda$$ есть ограниченная функция
на множестве $[0,\infty){\times}[0,\infty),$ то применение теоремы
Лебега о мажорируемой сходимости дает существование предела
функции $I_N(z)$ при $N\to\infty$ и выполнение равенства $$ V.P.
\int\limits_{-\infty}^\infty\frac{\varphi(\lambda)}{\lambda+z}\,d\lambda=
\int\limits_{-\infty}^\infty \left(\; V.P.
\int\limits_{-\infty}^\infty \frac{\cos \lambda t}{\lambda+z}
\,d\lambda \right)d\mu(t). $$

Тем самым теорема доказана.

Следующая теорема тесно связана с уже доказанными
теоремами~\ref{Skoryk-Theorem-4} и \ref{Skoryk-Theorem-5}. Её
можно считать одним из основных результатов парграфа. Напомним,
что определение преобразования Карлемана меры дано во введении.

\begin{theorem}\label{Skoryk-Theorem-6}\hskip-2mm.
Пусть мера $\mu\in M({\Bbb R}),$ $\widehat \mu$ --- её
преобразование Фурье, $z$
--- невещественное число. Тогда выполняется равенство
$$ \frac{i}{\sqrt{2\pi}}\;V.P.\int\limits_{-\infty}^\infty
\frac{\widehat \mu(\lambda)}{\lambda+z}\, d\lambda=F(z),
\n{Skoryk-formula-6}$$ где $F$ --- преобразование Карлемана меры
$\mu$.
\end{theorem}

{\sc Доказательство.} Пусть мера $\mu^\bigtriangleup$ определяется
равенством $\mu^\bigtriangleup(E) = \mu(-E).$ Обозначим
$$\mu_1=\frac{1}{2}(\mu+\mu^\bigtriangleup),\quad
\mu_2=\frac{1}{2}(\mu-\mu^\bigtriangleup).$$ Мера $\mu_1$
--- чётная, а $\mu_2$
--- нечётная, причём $\mu=\mu_1+\mu_2.$ Обе меры $\mu_1$ и
$\mu_2$ принадлежат алгебре $M(\Bbb R).$ Тогда
$$\widehat \mu(\lambda)=\widehat \mu_1(\lambda)+\widehat
\mu_2(\lambda)=\varphi_1(\lambda)+\varphi_2(\lambda).$$ Имеем
$$\varphi_1(\lambda)=\frac{1}{\sqrt{2\pi}}\int\limits_{-\infty}^\infty
e^{-i\lambda
t}d\mu_1(t)=\frac{1}{\sqrt{2\pi}}\int\limits_{-\infty}^\infty
\cos\lambda
t\,d\mu_1(t)=\frac{2}{\sqrt{2\pi}}\int\limits_0^{\;\;\infty
\;_\prime} \cos\lambda t\,d\mu_1(t),$$
$$\varphi_2(\lambda)=\frac{1}{\sqrt{2\pi}}\int\limits_{-\infty}^\infty
e^{-i\lambda
t}d\mu_2(t)=-\frac{i}{\sqrt{2\pi}}\int\limits_{-\infty}^\infty
\sin\lambda
t\,d\mu_2(t)=-\frac{2i}{\sqrt{2\pi}}\int\limits_0^\infty
\sin\lambda t\,d\mu_2(t).$$ Значение символа $^\prime$ над знаком
интеграла объяснялось во вступлении.

Из теорем~\ref{Skoryk-Theorem-5} и \ref{Skoryk-Theorem-4} следуют
равенства
$$V.P.\int\limits_{-\infty}^\infty\frac{\varphi_1(\lambda)}{\lambda+z}\,d\lambda=
\frac{2}{\sqrt{2\pi}}\int\limits_0^{\;\;\infty\;_\prime}
\left(V.P.\int\limits_{-\infty}^\infty\frac{\cos\lambda
t}{\lambda+z}\,d\lambda\right)d\mu_1(t),$$
$$\int\limits_{-\infty}^\infty\frac{\varphi_2(\lambda)}{\lambda+z}\,d\lambda=
-\frac{2i}{\sqrt{2\pi}}\int\limits_0^\infty
\left(\int\limits_{-\infty}^\infty\frac{\sin\lambda
t}{\lambda+z}\,d\lambda\right)d\mu_2(t).$$ Учитывая чётность меры
$\mu_1$ и нечётность меры $\mu_2,$ после сложения написанных
равенств получаем
$$V.P.\int\limits_{-\infty}^\infty\frac{\widehat
\mu(\lambda)}{\lambda+z}\,d\lambda=\frac{1}{\sqrt{2\pi}}\left(\int\limits_{-\infty}^\infty
\left(V.P.\int\limits_{-\infty}^\infty\frac{\cos\lambda
t}{\lambda+z}\,d\lambda\right)d\mu_1(t)+\right.$$
$$+\left.\int\limits_{-\infty}^\infty
\left(\int\limits_{-\infty}^\infty\frac{-i\sin\lambda
t}{\lambda+z}\,d\lambda\right)d\mu_2(t)\right)=\frac{1}{\sqrt{2\pi}}\left(\int\limits_{-\infty}^\infty
\left(V.P.\int\limits_{-\infty}^\infty\frac{e^{-i\lambda
t}}{\lambda+z}\,d\lambda\right)d\mu_1(t)\right.+$$
$$+\left.\int\limits_{-\infty}^\infty
\left(\int\limits_{-\infty}^\infty\frac{e^{-i\lambda
t}}{\lambda+z}\,d\lambda\right)d\mu_2(t)\right)=\frac{1}{\sqrt{2\pi}}\int\limits_{-\infty}^\infty
\left(V.P.\int\limits_{-\infty}^\infty\frac{e^{-i\lambda
t}}{\lambda+z}\,d\lambda\right)d\mu(t).\n{Skoryk-formula-7}
$$

Далее понадобится следующее утверждение.

\begin{theorem}\label{Skoryk-Theorem-7} \hskip-2mm. Пусть
$g(z)$ --- функция, мероморфная в полуплоскости $\Im z{>}0$ с
конечным числом полюсов в точках $a_1,$ $\ldots,$ $a_n.$ Пусть
функция $g$ непрерывно продолжается на границу полуплоскости и
удовлетворяет соотношению $g(z)=o(1)$ при $z\to \infty.$ Тогда
выполняется равенство
$$V.P.\int\limits_{-\infty}^\infty g(u)e^{iu}\,du=2\pi i
\sum\limits_{m=1}^n Res_{a_m} g(z)e^{iz}.
$$
\end{theorem}

Это задача 28.04 из \cite{Skoryk-Evgrafov}.

Продолжим доказательство теоремы~\ref{Skoryk-Theorem-6}.
Справедливо равенство

$$I=V.P.\int\limits_{-\infty}^\infty \frac{e^{-i\lambda
t}}{\lambda+z}\, d\lambda=\left\{\begin{array}{cl} \displaystyle
-{\rm sign}\;t \int\limits_{-\infty}^\infty \frac{e^{iu}}{u-tz}\;
du,& \mbox{если}\;\; t\ne0,\\-i\pi \,{\rm sign}\,y,&
\mbox{если}\;\;t=0.\end{array}\right.$$

Предположим, что $\Im z>0.$ Тогда по
теореме~\ref{Skoryk-Theorem-7} имеем
$$I=\left\{\begin{array}{cl} -2\pi
ie^{itz},& \mbox{если}\;\;t>0,\\
0,& \mbox{если}\;\;t<0.
\end{array}\right.$$

Если же  $\Im z<0,$ то
$$I=\left\{\begin{array}{cl}
0,& \mbox{если}\;\;t>0,\\
2\pi ie^{itz},& \mbox{если}\;\;t<0.
\end{array}\right.$$
Разбивая внешний интеграл из правой части равенства
(\ref{Skoryk-formula-7}) на интегралы по множествам $(-\infty,0),$
$\{0\},$ $(0, \infty)$ и подставляя найденное значение $I,$
получаем утверждение теоремы.

Далее в равенстве (\ref{Skoryk-formula-6}) мы хотим перейти к
пределу, когда $z=x+iy\to x.$ Для этого понадобится следующее
утверждение.

\begin{theorem}\label{Skoryk-Theorem-8} \hskip-2mm. Пусть $\cal L$ ---
ориентированная компактная гладкая жорданова кривая, возможно
замкнутая, $z$ --- произвольная внутренняя точка на кривой $\cal
L,$ $\nu$ --- единичная нормаль к кривой $\cal L$ в точке $z,$
расположенная слева от касательной, $z_\varepsilon=z\pm\varepsilon
\nu,$ ${\cal L}_\varepsilon={\cal L}\setminus l_\varepsilon,$
$l_\varepsilon$
--- связная компонента множества ${\cal L}\bigcap
C(z,\varepsilon),$ содержащая точку $z,$
$C(z,\varepsilon)=\{\zeta:\;|\zeta-z|<\varepsilon\},\;\varphi$
--- непрерывная функция на кривой $\cal L.$ Тогда
$$\lim\limits_{\varepsilon\to+0}\frac{1}{2\pi i}\left(\int\limits_{\cal L}\frac{\varphi(\zeta)}{\zeta-z_\varepsilon}\;d\zeta-
\int\limits_{{\cal
L}_\varepsilon}\frac{\varphi(\zeta)}{\zeta-z}\;d\zeta\right)=\mp\frac{1}{2}\;\varphi(z).$$

\end{theorem}

\begin{remark}\label{Skoryk-Remark_2}
Сформулированная теорема есть частный случай оригинального
результата Привалова, который в \cite[глава 5,
\S~1]{Skoryk-Privalov}  носит название основная лемма для
интегралов типа Коши. В оригинале рассматривается случай, когда
$\cal L$ --- спрямляемая кривая, $\varphi$ --- интегрируемая
функция на кривой $\cal L,$ $\nu$ --- единичный некасательный
вектор в точке $z,$ расположенный слева от касательной.
\end{remark}

\begin{theorem}\label{Skoryk-Theorem-9} \hskip-2mm. Пусть
$\mu\in M(\Bbb R)$ и $\widehat \mu$ --- её преобразование Фурье.
Тогда для любого вещественного $x$ существует интеграл
$$\int\limits_{-\infty}^\infty\frac{\widehat \mu(\lambda)}{\lambda-x}\,d\lambda,\n{Skoryk-formula-10}
$$
понимаемый как $$\lim\limits_{{\varepsilon\to
0}\atop{N\to\infty}}\left(\int\limits_{-N}^{x-\varepsilon}
\frac{\widehat
\mu(\lambda)}{\lambda-x}\,d\lambda+\int\limits_{x+\varepsilon}^N
\frac{\widehat \mu(\lambda)}{\lambda-x}\,d\lambda\right),$$ и
выполняется равенство
$$\frac{1}{\pi}\int\limits_{-\infty}^{\infty} \frac{\widehat
\mu(\lambda)}{\lambda-x}\,d\lambda=-i\widehat\nu(x),$$ где
$d\nu(t)=({\rm sign}\,t) d\mu(t).$
\end{theorem}

{\sc Доказательство.} Величина $$\int\limits_{N}^{\infty}\left(
\frac{\widehat \mu(\lambda)}{\lambda-x+i\varepsilon
}-\frac{\widehat \mu(\lambda)}{\lambda-x}\right)d\lambda$$ в силу
ограниченности функции $\widehat \mu(\lambda),$ очевидно,
стремится к нулю при $\varepsilon\to 0.$ Поэтому, если применять
теорему~\ref{Skoryk-Theorem-8} к функции
$\varphi(\lambda)=\widehat \mu(\lambda),$ то в качестве кривой
$\cal L$ можно брать вещественную ось. Тогда получаем
$$\lim\limits_{\varepsilon\to
0}\left(\int\limits_{-\infty}^\infty \frac{\widehat
\mu(\lambda)}{\lambda-x\pm
i\varepsilon}\,d\lambda-\int\limits_{{\cal L}_\varepsilon}
\frac{\widehat \mu(\lambda)}{\lambda-x}\,d\lambda\right)=
\left\{\begin{array}{r}
-\pi i \widehat \mu(x),\\
\pi i \widehat \mu(x).
\end{array}\right. \n{Skoryk-formula-11}$$
Подчеркнём, что в написанном выше равенстве оба интеграла следует
понимать в смысле главного значения. Из
теоремы~\ref{Skoryk-Theorem-6} следует, что
$$
\lim\limits_{\varepsilon\to 0} V.P.\int\limits_{-\infty}^\infty
\frac{\widehat \mu(\lambda)}{\lambda-x+
i\varepsilon}\,d\lambda=-i\sqrt{2\pi}\int\limits_0^{\;\;\;\infty\;\;\;_\prime}
e^{-ixt}\,d\mu(t),\n{Skoryk-formula-12}$$
$$
\lim\limits_{\varepsilon\to 0} V.P.\int\limits_{-\infty}^\infty
\frac{\widehat \mu(\lambda)}{\lambda-x-
i\varepsilon}\,d\lambda=i\sqrt{2\pi}\int\limits_{-\infty}^{\;\;\;0\;\;\;_\prime}
e^{-ixt}\,d\mu(t).\n{Skoryk-formula-13}$$

Обозначения со  штрихами объяснялись в конце первого параграфа.

 Из сказанного следует, что существует
$$\lim\limits_{\varepsilon\to 0} V.P.\int\limits_{{\cal L}_\varepsilon}
\frac{\widehat \mu(\lambda)}{\lambda-x}\,d\lambda,$$ т.е.
существует интеграл (\ref{Skoryk-formula-10}) в смысле, описанном
в условии теоремы.

Таким образом, каждое слагаемое в левой части
(\ref{Skoryk-formula-11}) имеет предел. Теперь из
(\ref{Skoryk-formula-11})--(\ref{Skoryk-formula-13}) следуют
равенства
$$- i\sqrt{2\pi}\int\limits_0^{\;\;\;\infty\;\;\;_\prime}
e^{-ixt}d\mu(t)-\int\limits_{-\infty}^\infty\frac{\widehat
\mu(\lambda)}{\lambda-x}\,d\lambda=-\pi i \widehat \mu(x),$$
$$ i\sqrt{2\pi}\int\limits_{-\infty}^{\;\;\;0\;\;\;_\prime}
e^{-ixt}d\mu(t)-\int\limits_{-\infty}^\infty\frac{\widehat
\mu(\lambda)}{\lambda-x}\,d\lambda=\pi i \widehat \mu(x).$$

Одним из следствий этих равенств является тривиальная формула
$$\widehat \mu(x)=\frac{1}{\sqrt{2\pi}}\int\limits_{-\infty}^\infty
e^{-itx}d\mu(t).$$ Для нас важно другое следствие, которое
совпадает с утверждением теоремы, и которое получается
суммированием написанных равенств.

Теорема доказана.

Алгебра $\KR$ становится банаховой алгеброй, если норму элемента
$\widehat\mu$ определять формулой
$\displaystyle\|\widehat\mu\|=\frac{1}{\sqrt{2\pi}}\|\mu\|.$ Как
банахово пространство алгебра $\KR$ является прямой суммой
одномерного пространства $A_1$ постоянных функций и пространства
$A_2$ тех функций из $\KR,$ которые являются преобразованиями
Фурье тех мер из $M(\Bbb R), $ которые не нагружают нуля.

Напомним, что для функций, определенных на вещественной оси,
преобразование Гильберта определяется формулой
$$(Hf)(x)=V.P.\frac{1}{\pi}\int
\limits_{-\infty}^{+\infty}\frac{f(t)}{x-t}dt.$$

Рассматривается и обобщенное преобразование Гильберта
$$(hf)(x)=V.P.\frac{1}{\pi}\int\limits_{-\infty}^{\infty}\left(\frac{1}{x-t}\;+\frac{t}{1+t^2}
\right)f(t)\;dt.$$

 Хорошо известно, что оператор Гильберта $H$ является
 изометрическим в пространстве $L_2(-\infty,\infty)$ и что для
 $f\in L_2(-\infty,\infty)$ выполняется равенство $H^2 f=-f.$
 Справедливо также следующее утверждение.

\begin{theorem}\label{Skoryk-Theorem-E1} \hskip-2mm.
Пусть измеримая на вещественной оси функция $f(x)$ такова, что
выполняется неравенство
$$
\int\limits_{-\infty}^\infty \frac{|f(x)|^2}{1+x^2}\;dx<\infty.
$$
 Тогда существует постоянная $C$ такая, что
$$
h^2 f = -f+C.
 $$
\end{theorem}

Доказательство теоремы~\ref{Skoryk-Theorem-E1} можно найти в
\cite[ пункт 3.6.2]{Skoryk-J}.

Опишем свойства оператора Гильберта в банаховом пространстве
$\KR,$ являющимся прямой суммой пространств $A_1$ и $A_2.$
Следующую теорему можно рассматривать как вариант
теоремы~\ref{Skoryk-Theorem-9}.

\begin{theorem}\label{Skoryk-Theorem-H1} \hskip-2mm.
На функциях из банахового пространства $\KR=A_1\oplus A_2$
преобразование Гильберта корректно определено для всех $x\in{\Bbb
R}.$ Оператор Гильберта $H$ переводит пространство $A_1$  в ноль.
На пространстве $A_2$ оператор $H$ является изометрическим, причём
$HA_2=A_2.$ Каждая функция из $A_2,$ которая является
преобразованием Фурье меры $\mu,$ сосредоточенной на одной из
полуосей $(-\infty,0)$ или $(0, \infty),$ является собственной
функцией оператора $H.$ Любая функция из $A_2$ является суммой
двух собственных функций оператора $H.$
\end{theorem}

Пусть $C_b$ --- пространство непрерывных ограниченных функций на
вещественной оси. Следующая теорема это критерий принадлежности
функции из пространства $C_b$ алгебре $\KR.$

\begin{theorem}\label{Skoryk-Theorem-H2} \hskip-2mm.
Для того, чтобы функция $f$ из пространства $C_b$ принадлежала
алгебре $\KR,$ необходимо и достаточно, чтобы для любой точки
$x\in{\Bbb R}$ преобразование Гильберта функции $f$ было корректно
определено и чтобы $Hf\in\KR.$
\end{theorem}

{\sc Доказательство.} Часть теоремы в сторону необходимости
следует из теоремы~\ref{Skoryk-Theorem-9}. Докажем вторую часть.
Пусть $\varphi=Hf\in\KR.$ По теореме~\ref{Skoryk-Theorem-H1}
существует функция $f_1\in\KR$ и постоянная $C_1$ такие, что
$Hf=\varphi=Hf_1+C_1.$ Следовательно, $H^2f=H^2f_1.$ Теперь из
теоремы \ref{Skoryk-Theorem-E1} следует, что  $f=f_1+C_2.$ Поэтому
$f\in\KR.$

Теорема доказана.

Ввиду важности алгебры $\FR$ сформулируем соответствующие
результаты для этой алгебры. Известно \cite[теорема
19.18]{Skoryk-Huit-1}, что алгебра $\FR$ является идеалом в
алгебре $\KR.$

\begin{theorem}\label{Skoryk-Theorem-H3} \hskip-2mm.
Пусть $f\in L_1(-\infty,\infty),$ $\widehat f$ --- её
преобразование Фурье, $z$ --- невещественное число. Тогда
выполняется равенство
$$\frac{i}{\sqrt{2\pi}}\;V.P.\int\limits_{-\infty}^\infty\frac{\widehat f(\lambda)}{\lambda+z}\;d\lambda
=F(z),$$ где $F(z)$ --- преобразование Карлемана функции $f.$
\end{theorem}

Это следствие теоремы~\ref{Skoryk-Theorem-6}.

\begin{theorem}\label{Skoryk-Theorem-H4} \hskip-2mm.
Пусть функция $f\in L_1(-\infty,\infty)$ и $\widehat f$ --- её
преобразование Фурье. Тогда для любого вещественного $x$
существует интеграл
$$ V.P.\int\limits_{-\infty}^\infty\frac{\widehat
f(\lambda)}{\lambda-x}\,d\lambda,
$$
понимаемый как $$\lim\limits_{{\varepsilon\to
0}\atop{N\to\infty}}\left(\int\limits_{-N}^{x-\varepsilon}
\frac{\widehat
f(\lambda)}{\lambda-x}\,d\lambda+\int\limits_{x+\varepsilon}^N
\frac{\widehat f(\lambda)}{\lambda-x}\,d\lambda\right),$$ и
выполняется равенство
$$V.P.\frac{1}{\pi}\int\limits_{-\infty}^{\infty} \frac{\widehat
f(\lambda)}{\lambda-x}\,d\lambda=-i\widehat f_1(x),$$ где
$f_1(t)=({\rm sign}\,t) f(t).$
\end{theorem}

Это следствие теоремы~\ref{Skoryk-Theorem-9}.

\begin{theorem}\label{Skoryk-Theorem-H5} \hskip-2mm.
Оператор Гильберта $H$ отображает банахово пространство $\FR$ на
$\FR.$ Это изометрический оператор. Всякая функция $\varphi$ из
$\FR,$ являющаяся преобразованием Фурье функции $f,$ равной нулю
на одной из полуосей $(-\infty,0),$ $(0, \infty),$ является
собственной функцией оператора $H.$ Любая функция из $\FR$
представляется в виде суммы двух собственных функций
оператора~$H.$
\end{theorem}

Это следствие теоремы~\ref{Skoryk-Theorem-H1}.

Пусть $C_0$ --- пространство непрерывных на оси $(-\infty,\infty)$
бесконечно малых на бесконечности функций.

\begin{theorem}\label{Skoryk-Theorem-H6} \hskip-2mm.
Пусть $\varphi\in C_0.$ Для того, чтобы функция $\varphi$
принадлежала алгебре $\FR,$ необходимо и достаточно, чтобы для
любого $x\in{\Bbb R}$ преобразование Гильберта функции $\varphi$
было корректно определено и функция $H\varphi$ принадлежала $\FR.$
\end{theorem}

Доказательство этой теоремы аналогично доказательству
теоремы~\ref{Skoryk-Theorem-H2}.

Рассмотрим еще аналоги
теорем~\ref{Skoryk-Theorem-H3}---\ref{Skoryk-Theorem-H6} для
банаховой алгебры Винера $\cal W,$ состоящей из функций вида
$$\varphi(x)=\sum\limits_{n=-\infty}^\infty c_n e^{inx},\quad \|\varphi\|=\sum\limits_{n=-\infty}^\infty|c_n|<\infty.$$

По сложившейся традиции символом $\varphi$  будем обозначать также
функцию $$\varphi(\zeta)=\sum\limits_{n=-\infty}^\infty c_n\zeta
^n,\quad \zeta=e^{ix}.$$

Это позволяет в одних случаях рассматривать $\varphi$ как
периодическую функцию на вещественной оси, а в других случаях как
функцию на единичной окружности ${\Bbb
T}=\{\zeta:\;\;|\zeta|=1\}.$ Отметим, что в дальнейшем там, где
окружность $\Bbb T$ рассматривается как ориентированная кривая,
считается, что она пробегается против хода часовой стрелки.

\begin{theorem}\label{Skoryk-Theorem-H7} \hskip-2mm.
Пусть
$$
\varphi(\zeta)=\sum\limits_{n=-\infty}^\infty c_n\zeta^n\;\in{\cal
W},\n{formula-17191}
$$
$$
\psi(z)=\int\limits_{\Bbb T}
\frac{\varphi(\zeta)}{\zeta-z}\;d\zeta.\n{formula-17192}
$$

Тогда
$$\psi(z)=\left\{\begin{array}{cl}
\displaystyle 2\pi i\sum\limits_{n=0}^\infty c_n z^n,
&\mbox{если}\;\;|z|<1,\\[10pt]
\displaystyle -2\pi i\sum\limits_{n=-\infty}^{-1} c_n z^n,
&\mbox{если}\;\;|z|>1.
\end{array}\right.$$

\end{theorem}

{\sc Доказательство.} Если в формулу (\ref{formula-17192})
подставить $\varphi(\zeta)$ из формулы (\ref{formula-17191}), то в
силу равномерной сходимости ряда порядок интегрирования и
суммирования можно поменять. Отсюда уже легко следует утверждение
теоремы.~

\begin{remark} \hskip-2mm.
Если стать на точку зрения абстрактного гармонического анализа, то
функцию $\varphi$ можно считать преобразованием Фурье функции
$c_n,$ определённой на абелевой группе $\Bbb Z$ целых чисел. Тогда
функцию $\psi$ можно рассматривать как преобразование Карлемана
функции $c_n.$ Таким образом, теорема~\ref{Skoryk-Theorem-H7}
является  аналогом теоремы~\ref{Skoryk-Theorem-H3} для алгебры
$\cal W.$
\end{remark}

\begin{theorem}\label{Skoryk-Theorem-H8} \hskip-2mm.
Пусть $\displaystyle \varphi=\sum\limits_{n=-\infty}^\infty
c_ne^{inx}\in{\cal W}.$ Тогда для любого $z\in{\Bbb T}$
выполняется равенство
$$
V.P.\frac{1}{\pi}\int\limits_{\Bbb T}
\frac{\varphi(\zeta)}{\zeta-z}\;d\zeta=i\sum\limits_{n=-\infty}^{\infty}\eta_n
c_n z^n,\n{formula-05061}
$$
где $\eta_n=1$ при $n\ge0$ и $\eta_n=-1$ при $n<0.$
\end{theorem}

{\sc Доказательство.} Согласно теореме~\ref{Skoryk-Theorem-8}
выполняется равенство
$$\lim\limits_{\varepsilon\to0}\frac{1}{2\pi i}\left(\int\limits_{\Bbb T}\frac{\varphi(\zeta)}{\zeta-z_\varepsilon}
\;d\zeta-\int\limits_{{\Bbb
T}_\varepsilon}\frac{\varphi(\zeta)}{\zeta-z}\;d\zeta\right)=\mp\frac{1}{2}\varphi(z).$$

Из теоремы~\ref{Skoryk-Theorem-H7} следует, что
$$\lim\limits_{\varepsilon\to+0}\frac{1}{2\pi i}\int\limits_{\Bbb T}\frac{\varphi(\zeta)}{\zeta-z+\varepsilon e^{i\theta}}\;d\zeta
=\sum\limits_{n=0}^\infty c_n z^n,\quad\theta=\arg z,$$
$$\lim\limits_{\varepsilon\to+0}\frac{1}{2\pi i}\int\limits_{\Bbb T}\frac{\varphi(\zeta)}{\zeta-z-\varepsilon e^{i\theta}}\;d\zeta
=-\sum\limits_{n=-\infty}^{-1} c_n z^n.$$

Из сформулированных утверждений следует существование интеграла
(\ref{formula-05061}) и выполнение равенств
$$2\pi i\sum_{n=0}^\infty c_n z^n -V.P.\int\limits_{\Bbb T}\frac{\varphi(\zeta)}{\zeta-z}\;d\zeta=\pi i\varphi(z),$$
$$-2\pi i\sum_{n=-\infty}^{-1} c_n z^n -V.P.\int\limits_{\Bbb T}\frac{\varphi(\zeta)}{\zeta-z}\;d\zeta=-\pi
i\varphi(z).$$ Из этого следует утверждение теоремы.

Теорема доказана.

Пусть ${\cal W}_0$ --- подпространство банахового пространства
$\cal W,$ состоящее из тех функций $w\in{\cal W},$ для которых
$c_0=0.$ В пространстве ${\cal W}_0$ оператор $\displaystyle
w\longmapsto V.P.\frac{1}{\pi}\int\limits_{\Bbb
T}\frac{w(\zeta)}{z-\zeta}\;d\zeta$ совпадает с оператором
Гильберта $H,$ который согласно общепринятому определению имеет
вид
$$(Hw)(x)=V.P.\frac{1}{2\pi}\int\limits_0^{2\pi}\ctg\frac{x-y}{2}\;w(y)dy.$$

\begin{theorem}\label{Skoryk-Theorem-H9}\hskip-2mm.
Оператор Гильберта $H$ отображает банахово пространство ${\cal
W}_0$ на себя. В пространстве ${\cal W}_0$ это изометрический
оператор. Всякая функция $w$ из ${\cal W}_0$ такая, что $c_n=0$
либо при $n>0,$ либо при $n<0$ является собственной функцией
оператора $H.$ Любая функция из ${\cal W}_0$ представляется в виде
суммы двух собственных функций оператора $H.$
\end{theorem}

Сформулированная теорема является лёгким следствием предыдущей.

\begin{theorem}\label{Skoryk-Theorem-H10} \hskip-2mm.
Пусть $w$ --- непрерывная $2\pi$-периодическая функция. Для того,
чтобы $w$ принадлежала алгебре $\cal W,$ необходимо и достаточно,
чтобы преобразование Гильберта этой функции было корректно
определено в каждой точке $x$ вещественной оси и чтобы $Hw\in{\cal
W}.$
\end{theorem}

{\sc Доказательство.} Часть теоремы в сторону необходимости
следует из теоремы~\ref{Skoryk-Theorem-H8}. Докажем вторую часть.
Пусть $Hw\in{\cal W}.$ По теореме~\ref{Skoryk-Theorem-H9}
существуют функция $w_1\in{\cal W}_0$ и постоянная $c$ такие, что
$H(w-w_1)=c.$ Функция $w-w_1\in L_2(0,2\pi).$ Для функций $u\in
L_2(0,2\pi)$ коэффициенты Фурье $c_n$ функции $u$ и коэффициенты
Фурье $d_n$ функции $Hu$ связаны соотношением $\displaystyle
\sum\limits_{n=-\infty \atop{n\neq0}}^\infty
|c_n|^2=\sum\limits_{n=-\infty}^\infty |d_n|^2,$ кроме того
$d_0=0.$ Отсюда следует, что $c=0$ и что все коэффициенты Фурье
функции $w-w_1$ кроме нулевого равны нулю. Поэтому функция $w-w_1$
является постоянной, $w\in{\cal W}.$

Теорема доказана.

\begin{remark} \hskip-2mm.
Отметим, что из теоремы~\ref{Skoryk-Theorem-H10} следует, что если
$w\in{\cal W},$ то $w$ и $Hw$ являются непрерывными функциями. Как
следует из утверждения~14 из вступления из непрерывности функций
$w$ и $Hw$ нельзя заключить, что $w\in{\cal W}.$

Теорема~\ref{Skoryk-Theorem-H10} по существу эквивалентна
утверждению из \cite[гл.~2, пункт~10]{Skoryk-Kahan},  где
говорится, что если функция $w$ принадлежит алгебре $\cal W$
(алгебре $A$ в обозначениях Кахана), то её сопряжённая функция
также принадлежит $\cal W.$
\end{remark}

Сейчас мы докажем аналог теоремы~\ref{Skoryk-Theorem-H3} для
функций из $L_p(\Bbb R)$.

\begin{theorem}\label{Skoryk-Theorem-11} \hskip-2mm.
Пусть функция  $f\in L_p(\Bbb R), \; p\in(1,2]$, $\widehat f$ ---
ее преобразование Фурье, $z$ --- невещественное число. Тогда
выполняется равенство
$$
\frac{i}{\sqrt{2\pi}}\int\limits_{-\infty}^{+\infty}\frac{\widehat
f(\lambda)}{\lambda+z}d\lambda=F(z),\n{virazhenie13.1}
$$
где $F(z)$ --- преобразование Карлемана функции $f$.
\end{theorem}

{\sc Доказательство.} Во-первых, заметим, что интеграл в левой
части равенства корректно определен. Действительно, пусть число
$q$ определяется из равенства
$\displaystyle\frac{1}{p}+\frac{1}{q}=1$. Тогда по теореме
Хаусдорфа-Юнга  имеем $\widehat f(\lambda)\in L_q.$ Кроме того,
$\displaystyle\frac{1}{\lambda+z}\in L_p$. Поэтому
$\displaystyle\frac{\widehat f(\lambda)}{\lambda+z}\in L_1$. Пусть
$\{f_n\}_{n=1}^\infty\subset  L_p\bigcap L_1$ такая
последовательность, что $f_n\to f$ в метрике $L_p$, пусть $F_n(z)$
--- преобразование Карлемана функции $f_n.$ По теореме~\ref{Skoryk-Theorem-H3}
$$\displaystyle\frac{i}{\sqrt{2\pi}}\int\limits_{-\infty}^{+\infty}\frac{\widehat
f_n(\lambda)}{\lambda+z}d\lambda=F_n(z). \n{virazhenie13.2}
$$

Пусть $\Im z>0.$ Имеем
$$\displaystyle|F(z)-F_n(z)|=\left|\int\limits_0^\infty(f(t)-f_n(t))e^{itz}dt\right|\le\|f-f_n\|_p\|e^{itz}\|_q=\left(\frac{1}{qy}\right)^\frac{1}{q}\|f-f_n\|_p.$$
Следовательно, $F_n(z)\to F(z) $ ($n\to\infty$). Кроме того,
 $$\displaystyle\left|\int\limits_{-\infty}^{+\infty}\frac{\widehat f(\lambda)}{\lambda+z}d\lambda-\int\limits_{-\infty}^{+\infty}\frac{\widehat f_n(\lambda)}{\lambda+z}d\lambda\right|\le\|\widehat f-
 \widehat f_n\|_q\left\|\frac{1}{\lambda+z}\right\|_p.$$
 Из этого следует, что
 $$\displaystyle\lim\limits_{n\to\infty}\int\limits_{-\infty}^{+\infty}\frac{\widehat f_n(\lambda)}{\lambda+z}d\lambda=
 \int\limits_{-\infty}^{+\infty}\frac{\widehat f(\lambda)}{\lambda+z}d\lambda.$$
 Переходя в равенстве (\ref{virazhenie13.2}) к пределу при
 $n\to\infty$, получим (\ref{virazhenie13.1}). Мы рассмотрели случай
$\Im z>0.$ Случай $\Im z<0$ рассматривается аналогично.

Теорема доказана.

\sect{Формула Повзнера}\label{4} В этом параграфе мы даем
интегральное представление для преобразования Карлемана мер из
классов $M_k$. Как уже отмечалось во вступлении, класс мер $M_k$
аналогичен классу функций $F_k,$ который ввёл Бохнер
\cite{Skoryk-Bohner}.
 Для мер $\mu\in M_k$ $k$-тое преобразование
Бохнера $\widehat \mu(k,t)$ определяется по формуле
$$\widehat \mu(k,t)=\frac{1}{\sqrt{2\pi} }\int\limits_{-\infty}^\infty
\frac{e^{-ixt}-P_{k-1}(x, t)}{(-ix)^k}\,d\mu(x),$$ где
$$P_{k-1}(x, t)=\left\{\begin{array}{cl}
\displaystyle
\sum\limits_{m=0}^{k-1}\frac{(-ixt)^m}{m!},&\mbox{если}\;\;|x|\le
1,\\[5pt]
0,& \mbox{если}\;\;|x|> 1.
\end{array}\right.$$

Заметим, что $\widehat \mu(0,t)$ совпадает с преобразованием Фурье
$\widehat \mu(t)$ меры $\mu$. Если $\mu\in M_0,$ то $\mu\in M_k$
для любого $k.$ В этом случае выполняется равенство
$\;\displaystyle\frac{d^k\widehat \mu(k,t)}{dt^k}=\widehat
\mu(t).\;$ Отметим ещё следующее. Из известной формулы
$$f(x)=\sum\limits_{m=0}^{k-1}\frac{f^{(m)}(0)}{m!}\;x^m+\frac{1}{(k-1)!}\int
\limits_0^x f^{(k)}(u)(x-u)^{k-1}du,$$ применённой к функции
$f(x)=e^{itx},$ следует оценка
$$\left|\frac{1}{(-ix)^k}\left(e^{-itx}-\sum\limits_{m=0}^{k-1}\frac{(-itx)^m}{m!}\right)\right|\le\frac{|t|^k}{k!}.$$
Из этого неравенства, в свою очередь, легко следует, что если мера
$\mu$ из класса $M_k$ не нагружает нуля, то выполняется равенство
$$\lim\limits_{t\to\pm\infty}\frac{\widehat\mu(k,\; t)}{t^k}=0.$$

\begin{theorem}\label{Skoryk-Theorem-12} \hskip-2mm.
Пусть $\mu\in M_k,$ $\widehat \mu(k, \lambda)$ --- ее $k$-ое
преобразование Бохнера, $F(z)$ --- ее преобразование Карлемана,
$z$ --- невещественное число. Тогда справедлива формула
$$ \displaystyle\frac{i
k!}{\sqrt{2\pi}}V.P.\int\limits_{-\infty}^{+\infty}\frac{\widehat
\mu(k, \lambda)}{(\lambda+z)^{k+1}}d\lambda=F(z).\n{virazhenie14}
$$
\end{theorem}

{\sc Доказательство.} Для $k=0$ теорема~\ref{Skoryk-Theorem-12}
совпадает с теоремой~\ref{Skoryk-Theorem-6}. Поэтому в дальнейшем
можно считать, что $k\ge 1$. Заметим, что если формула
(\ref{virazhenie14}) верна для каждой из мер $\mu_1$ и $\mu_2$, то
она верна и для их суммы. Если $\mu$ --- это мера Дирака $\delta,$
то формула (\ref{virazhenie14}) доказывается прямым вычислением
левой и правой частей. Из сказанного следует, что формулу
(\ref{virazhenie14}) достаточно доказывать для случая, когда мера
$\mu$ не нагружает нуля. В дальнейшем считается, что это условие
выполняется.

Предположим, что $\mu\in M_0$. Тогда $\mu\in M_k$ для любого $k$.
В этом случае формула (\ref{virazhenie14}) получается из формулы
(\ref{Skoryk-formula-6}) с помощью интегрирования по частям.

 Пусть $\mu$
--- произвольная мера из $M_k$. Тогда $\mu=\mu_1+\mu_2$, где
$\mu_1$ --- ограничение $\mu$ на сегмент $[-1, 1],$ а
$\mu_2=\mu-\mu_1.$ Каждая из мер $\mu_1,$ $\mu_2$ принадлежит
$M_k$, причём $\mu_1\in M_0$. По доказанному, формула
(\ref{virazhenie14}) верна для меры $\mu_1$.

Из сказанного следует, что для справедливости формулы
(\ref{virazhenie14}) достаточно проверить ее справедливость при
выполнении двух дополнительных условий:

1) мера $\mu$ не нагружает сегмент $[-1, 1]$,

2) $k\ge 1$.

Далее мы считаем, что эти условия выполняются. В этом случае
$$\displaystyle\int\limits_{-\infty}^{+\infty}\frac{\widehat \mu(k, \lambda)}{(\lambda+z)^{k+1}}d\lambda
=\frac{1}{\sqrt{2\pi}}\int\limits_{-\infty}^{+\infty}\left(\int\limits_{-\infty}^{+\infty}
\frac{e^{-it\lambda}}{(-it)^k(\lambda+z)^{k+1}}d\mu(t)\right)d\lambda.$$
Выполняется неравенство
$$\displaystyle\int\limits_{-\infty}^{+\infty}\left(\int\limits_{-\infty}^{+\infty}\left|
\frac{e^{-it\lambda}}{(-it)^k(\lambda+z)^{k+1}}\right|d|\mu|(t)\right)d\lambda<\infty$$
Поэтому из теорем Тонелли и Фубини следует, что
$$
\displaystyle\int\limits_{-\infty}^{+\infty}\frac{\widehat \mu(k,
\lambda)}{(\lambda+z)^{k+1}}d\lambda=
\frac{1}{\sqrt{2\pi}}\int\limits_{-\infty}^{+\infty}\left(\frac{1}
{(-it)^k}\int\limits_{-\infty}^{+\infty}\frac{e^{-it\lambda}}{(\lambda+z)^{k+1}}
d\lambda\right)d\mu(t). \n{virazhenie15} $$ Делая замену
$-t\lambda=u$, получим
$$\displaystyle\int\limits_{-\infty}^{+\infty}\frac{e^{-it\lambda}}{(\lambda+z)^{k+1}}d\lambda=
-sign \,t \;(-t)^k
\int\limits_{-\infty}^{+\infty}\frac{e^{iu}}{(u-tz)^{k+1}}du.$$

Предположим, что $\Im z>0.$ Тогда по
теореме~\ref{Skoryk-Theorem-7} имеем
$$\displaystyle\int\limits_{-\infty}^{+\infty}\frac{e^{iu}}{(u-tz)^{k+1}}du=0, \quad t<0.$$

Если же $t>0,$ то
$$\int\limits_{-\infty}^{+\infty}\frac{e^{iu}}{(u-tz)^{k+1}}du=
2\pi i{\rm Res}_{tz} \frac{e^{iu}}{(u-tz)^{k+1}}=$$ $$=2\pi i
e^{itz}{\rm Res}_{tz}\frac{e^{i(u-tz)}}{(u-tz)^{(k+1)}} =2\pi i
e^{itz}\frac{i^k}{k!}.$$

Подставляя в формулу (\ref{virazhenie15}) вычисленное значение
интеграла, получим
$$\displaystyle\int\limits_{-\infty}^{+\infty}\frac{\widehat \mu(k, \lambda)}{(\lambda+z)^{k+1}}d\lambda
=-\frac{1}{\sqrt{2\pi}}\int\limits_0^\infty\frac{1}{(-it)^k}2\pi
ie^{itz}(-it)^k\frac{1}{k!}d\mu(t)=-\frac{2\pi
i}{k!\sqrt{2\pi}}F(z).$$

Тем самым формула (\ref{virazhenie14}) доказана при $\Im z>0.$ В
случае $\Im z<0$ доказательство проводится аналогично.

\begin{remark} \hskip-2mm.
В \cite{Skoryk-Povzner} формула (\ref{virazhenie14}) приводится
для случая, когда мера $\mu$ имеет вид $f(t)dt,$ где $f(t)$ ---
непрерывная ограниченная функция (следовательно, $\mu\in M_2$).
Поэтому формулу (\ref{virazhenie14}) мы называем формулой
Повзнера.
\end{remark}

Сформулируем еще вариант теоремы~\ref{Skoryk-Theorem-12} для
классов Бохнера $F_k.$

\begin{theorem}\label{Skoryk-Theorem-12-1} \hskip-2mm.
Пусть функция $f\in F_k,$ $ \widehat f(k, \lambda)$ --- ее $k$-ое
преобразование Бохнера, $F(z)$ --- ее преобразование Карлемана,
$z$ --- невещественное число. Тогда справедлива формула
$$\frac{i
k!}{\sqrt{2\pi}}V.P.\int\limits_{-\infty}^{+\infty}\frac{\widehat
f(k, \lambda)}{(\lambda+z)^{k+1}}d\lambda=F(z).$$
\end{theorem}

\vskip5mm
\begin{large}
\noindent {{\bf A.\,F.~Grishin (А.\,Ф.~Гришин)}\\
V.~N.~Karazin Kharkov National University,\\
4 Svobody sq., Kharkov 61077, Ukraine}\\
{\it E-mail:} {grishin@univer.kharkov.ua}\\
\\
\noindent {{\bf M.\,V.~Skoryk (М.\,В.~Скорик)}\\
V.~N.~Karazin Kharkov National University,\\
4 Svobody sq., Kharkov 61077, Ukraine}\\
{\it E-mail:} {maksym\_skoryk@mail.ru}
\end{large}

\end{Large}
\end{document}